\documentclass{article}
\usepackage{fullpage}
\usepackage{graphicx}
\usepackage{amsmath}
\usepackage{amsfonts}
\usepackage{bbm}
\usepackage{amssymb}
\usepackage{amsthm}
\usepackage{subfigure}
\usepackage{cite}
\newtheorem{thm}{Theorem}
\newtheorem{deff}{Definition}

\begin{document}
\title{Numerical stability of explicit Runge-Kutta finite-difference schemes for the nonlinear Schr{\"o}dinger equation}
\author{R. M. Caplan\footnote{Corresponding author.  Present address:  Predictive Science Inc.  9990 Mesa Rim Rd, Suite 170, San Diego, CA 92121. email: caplanr@predsci.com, phone: 858-225-2314}~~and R. Carretero-Gonz{\'a}lez\\[1.0ex]
Nonlinear Dynamical System Group\footnote{\texttt{URL}: http://nlds.sdsu.edu},
Computational Science Research Center, and\\
Department of Mathematics and Statistics,
San Diego State University,\\
San Diego, California 92182-7720, USA\\
}
\date{\today}
\maketitle
\begin{abstract}
Linearized numerical stability bounds for solving the nonlinear time-dependent Schr{\"o}dinger equation (NLSE) using explicit finite-differencing are shown. The bounds are computed for the fourth-order Runge-Kutta scheme in time and both second-order and fourth-order central differencing in space.  Results are given for Dirichlet, modulus-squared Dirichlet, Laplacian-zero, and periodic boundary conditions for one, two, and three dimensions.  Our approach is to use standard Runge-Kutta linear stability theory, treating the nonlinearity of the NLSE as a constant.  The required bounds on the eigenvalues of the scheme matrices are found analytically when possible, and otherwise estimated using the Gershgorin circle theorem.
\end{abstract}

\section{Introduction}
The nonlinear Schr{\"o}dinger equation (NLSE) is used to model a wide variety of physical systems since it describes, to least nonlinear order, modulated wave propagation \cite{NLSE_nlpdebook}.  The general form of the NLSE can be written as
\begin{equation}
\label{nlse}
i\frac{\partial \Psi}{\partial t} + a\nabla^2\Psi - V({\bf r})\Psi + s|\Psi|^2 \Psi = 0,
\end{equation}
where $\Psi \in \mathbbm{C} $ is the value of the wavefunction, $\nabla^2$ is the Laplacian operator, and where $a>0$ and $s$ are parameters defined by the system being modeled.  $V(\bf{r})$ is an external potential term, which when included, makes Eq.~(\ref{nlse}) known as the Gross-Pitaevskii equation \cite{BEC_RCbook}.

Often, efficient and easy-to-use numerical methods are employed to simulate the NLSE. One such method is the method of lines where the time-stepping and spatial differencing are treated independently. This transforms the partial differential equation (PDE) into a large number of coupled ordinary differential equations (ODEs). These ODEs can then be solved using a variety of numerical schemes, one of the most common being the fourth order Runge-Kutta (RK4) scheme \cite{RK4}. Using the RK4 scheme with the NLSE produces a fully explicit scheme where each grid point at time $t$ is only a function of values at time $t-{\bf k}$ where ${\bf k}$ is the time-step. This simplifies computational implementations because no matrices are needed to be formed and stored, and no linear systems are needed to be solved (which in the nonlinear case also require a nonlinear iterative process). 

The only drawback to using explicit finite-difference schemes (such as the RK4) for simulating PDEs is that they are conditionally stable.  This means that there is an upper bound on the allowed size of the time-step which is dependent on the spatial-step size. If the time-step is larger than this bound, the scheme is unstable and diverges \cite{FD_PDE_BOOK}. Although rough estimates of the stability bound can be found through an inefficient educated guess-and-check, for higher dimensional scenarios, as well as long and/or large simulations, a more refined and predictable stability bound is essential for efficient simulations.   

In this paper, we formulate linearized stability bounds for simulating the NLSE with the RK4 scheme. The stability bounds depend on the specific spatial differencing scheme being used, as well as on the boundary conditions. We formulate the bounds for both second-order and fourth-order spatial differencing with a variety of boundary conditions (Dirichlet, modulus-squared Dirichlet (MSD), Laplacian-zero (L0), and periodic). Each analysis is done for one, two and three dimensions.

The paper is organized as follows. In Sec.~\ref{s:rk4sb} we review the basic RK4 stability properties and apply the results to the NLSE to formulate general stability bounds. Our basic procedure in finding the stability bounds and the linear algebra theorems that we utilize are also discussed. In Sec.~\ref{s:bc} we summarize the forms of the boundary conditions we consider. Our main analysis begins in Sec.~\ref{s:1dstb} with the one-dimensional NLSE. Linearized stability bounds are found for each scheme and boundary condition combination. In Sec.~\ref{s:2d} and \ref{s:3dstb}, we use the same procedures to formulate the bounds for the two- and three-dimensional NLSE respectively.  A few numerical examples showing the accuracy of the bounds are shown in Sec.~\ref{s:num}.  In Sec.~\ref{s:sum}, we conclude and summarize all the results from Secs.~\ref{s:1dstb}, \ref{s:2d}, and \ref{s:3dstb} into a concise reference.

\section{Stability theory}
\label{s:rk4sb}
\subsection{General Runge-Kutta scheme stability}
\label{s:rk4sb_gen}
Given an initial value problem of a set of linear first-order ODEs (in our case, a method-of-lines explicit PDE finite-difference scheme), one can formulate the matrix notation
\begin{equation}
\label{scheme}
\frac{\partial \vec \Psi}{\partial t} = {\cal A}\, \vec \Psi,
\end{equation}
where ${\cal A}$ contains the coefficients of the right-hand-sides of the ODEs.  We now define
\begin{equation}
\label{hhatdef}
\vec p = {\bf k}\,\vec \lambda,
\end{equation}
where ${\bf k}$ is the time-step size and $\vec \lambda$ contains the eigenvalues of ${\cal A}$. In our case, the eigenvalues of ${\cal A}$ will have the spatial-step size (denoted $h$) included in them, as well as any parameters of the NLSE.  As shown in Ref.~\cite{RK4_STB_BOOK}, for the fourth-order Runge-Kutta scheme, if a vector $\vec R(\vec p)$ is defined whose elements are the polynomials
\begin{equation}
\label{req}
R(p) = 1 + p + \frac{p^2}{2} + \frac{p^3}{6} + \frac{p^4}{24},
\end{equation}
then the stability of the RK4 scheme is guaranteed if
\begin{equation}
\label{rbound}
\lVert \vec R(\vec p)\rVert_\infty < 1,
\end{equation}
where $\lVert \; \rVert_{\infty}$ denotes the infinity norm defined as $\lVert \vec x \rVert_{\infty} = \mbox{max}\{\lvert x_0 \rvert, \lvert x_1 \rvert, ... ,\lvert x_{N-1} \rvert\}$.  Inserting Eq.~(\ref{hhatdef}) into Eq.~(\ref{req}) yields 
\begin{alignat}{2}
\label{req2}
\left|R(\lambda)\right|^2 = 1 & + \frac{1}{576}\,{\bf k}^8 |\lambda|^8 -
 \frac{1}{72}\,{\bf k}^6 |\lambda|^6 + \left(\frac{{\bf k}^6|\lambda|^6}{6} - {\bf k}^4| \lambda|^4 + 24 \right)\,\frac{{\bf k}}{12}\,\mbox{Re}(\lambda)    \\
&+ \left({\bf k}^4|\lambda|^4 + 24\right)\,\frac{{\bf k}^2}{12}\,(\mbox{Re}(\lambda))^2 + \left({\bf k}^2|\lambda|^2 + 4\right)\,\frac{{\bf k}^3}{3}\,(\mbox{Re}(\lambda))^3 + \frac{2{\bf k}^4}{3}\,(\mbox{Re}(\lambda))^4.
\notag
\end{alignat}

\begin{figure}[htbp]
\centering
\includegraphics[width=2.25in]{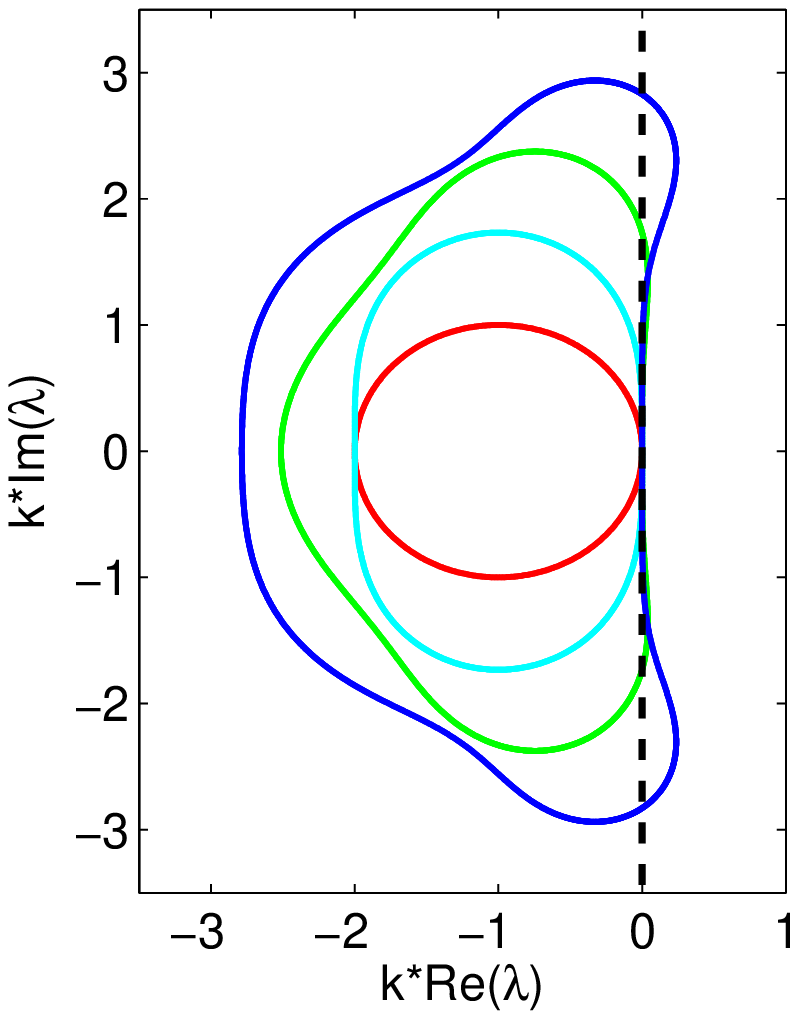}
\includegraphics[width=3.5in]{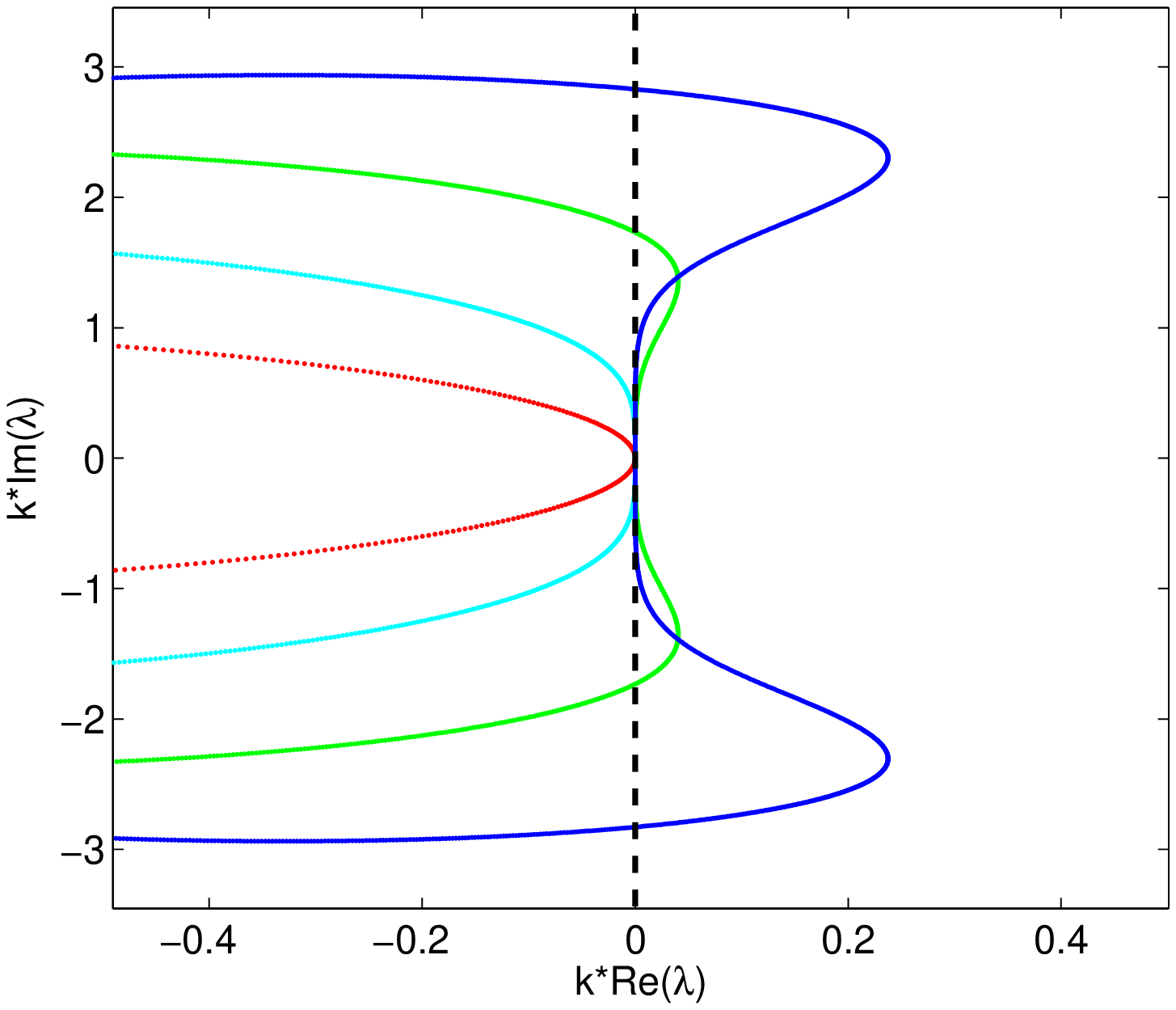}
\caption[Runge-Kutta stability regions.]{(Color online) Left:  Stability regions for Runge-Kutta schemes. Schemes from first-order to fourth-order are shown from center outwards. Right:  Magnified view of the same plot near the point where $\mbox{Re}(\lambda) = 0$.\label{f:rk4stb}}
\end{figure}

In Fig.~\ref{f:rk4stb} we show the stability region for the RK4 scheme given by Eq.~(\ref{rbound}) as well as that for lower-order Runge-Kutta schemes (whose $R(p)$ is defined by progressively truncated versions of Eq.~(\ref{req})) \cite{RK4_STB_BOOK}.
As we shall show, the eigenvalues of the $\cal A$ matrix are all purely imaginary (or nearly so) in the case of the nonlinear (and linear) Schr{\"o}dinger equation. Thus, it can be seen from Fig.~\ref{f:rk4stb} that the third-order Runge-Kutta is the lowest-order RK scheme that is conditionally stable for the Schr{\"o}dinger equations (however, as shown in Ref.~\cite{FD_NLSEGP_EXP_STB}, this is not the case if the real and imaginary parts of the NLSE are computed in a staggered time grid, or, as in Ref.~\cite{FD_STBNLSEDISS}, an artificial dissipative term is added to the NLSE), with the RK4 yielding a significantly larger bound on ${\bf k}$.  This is in contrast to similar PDEs such as the heat equation, whose $\cal A$ matrix eigenvalues are typically all real-valued, in which case even forward differencing (RK1) is conditionally stable.

If, as in our case, $\mbox{Re}(\vec \lambda)=\vec 0$, Eq.~(\ref{req2}) simplifies greatly and becomes
\begin{equation}
\label{req3}
\left|R(\vec \lambda)\right|^2 = 1 + \frac {1}{576}\,{\bf k}^8|\vec \lambda|^8-\frac{1}{72}\,{\bf k}^6|\vec \lambda|^6,
\end{equation}
in which case, Eq.~(\ref{rbound}) leads to the simple stability bound
\begin{equation}
\label{rk4stbeasy}
{\bf k} < \frac{\sqrt{8}}{\lVert\vec \lambda\rVert_\infty}.
\end{equation}
\subsection{Application to the NLSE}
Applying the above stability theory to the NLSE has the obvious problem that the analysis is purely linear, while the NLSE has one (or more) nonlinear terms.  A full nonlinear stability analysis is beyond the scope of this paper, so instead we linearize the problem by treating the nonlinearity ($|\Psi|^2$) as a constant value $U$ (the external potential term is usually a constant independent of $\Psi$ at each grid point, and so does not need any special treatment).  This has been done previously for the one-dimensional coupled NLSE for fourth-order differencing (in the exclusive case where $s<0$) in Ref.~\cite{RK4_2CNLSE_STB}.  Since the value of $|\Psi|^2$ changes over time during the simulation, the linearized stability bound will also change over time.  This change in many cases is expected to be small (which we have confirmed in numerical simulations, not reported here) and therefore can be ignored, i.e. one may compute the bound using the initial condition of $\Psi$ (and $V(\bf{r})$) and just leave a few percent leeway to cover any changes.  This is especially true in the repulsive case ($s<0$) where most situations have a constant-density background (or maximum background) and the dynamics do not cause the maximum background value to change significantly (for example, when simulated coherent structures, most of the dynamics are translations of the initial condition with little change in structure).  In attractive cases ($s>0$), blow-up can occur which can alter the stability bound greatly, causing the simulation to crash (although in such a case the wavefunction is exploding towards infinity, which most finite-difference schemes cannot handle anyways). Many times, simulations of a steady-state or near-steady-state in the modulus-squared with a constant potential are performed.  In such situations, the linearized stability bounds will be (nearly) exact. 

It is also useful to formulate stability bounds for the \emph{linear} Schr{\"o}dinger equation (LSE) (where $s=0$ and $V({\bf r})=0$).  In addition to providing bounds for the LSE, as will be discussed below, the results can also be used as practical estimates of the stability bounds for the NLSE (the discrepancy can often be solved by lowering the bound by a few percent).

\subsection{Stability analysis procedure}
In order to simplify the analysis, we first rewrite Eq.~(\ref{scheme}) as
\[
\frac{\partial \vec \Psi}{\partial t} = {\cal A} \vec \Psi = \frac{i\,a}{h^2}\,A \vec \Psi,
\]
where $h$ is the step-size of the spatial finite-difference scheme being used.  Then, assuming all eigenvalues of $A$ are real-valued, the stability condition of Eq.~(\ref{rk4stbeasy}) becomes
\begin{equation}
\label{rk4stbeasy2}
{\bf k} < \frac{\sqrt{8}}{\lVert\vec \lambda_{A}\rVert_\infty}\frac{h^2}{a}.
\end{equation}

In order to be able to use the stability bound of Eq.~(\ref{rk4stbeasy2}), we must first confirm that all eigenvalues of $A$ are purely real (or nearly so) for each scheme/boundary condition combination.  In cases where the eigenvalues are not able to be easily computed analytically, we show that the $A$ matrix's eigenvalues are a set of boundary values with the remaining eigenvalues being those of a symmetric matrix denoted $A^{\prime}$.  Then by  Thm.~\ref{t:real}, it is known that all the eigenvalues of $A$ are real.
\begin{thm}[Ref.~\cite{MATRIX_OP_BOOK}]
\label{t:real}
The eigenvalues of a real symmetric matrix are real.
\end{thm}

Once it has been established that Eq.~(\ref{rk4stbeasy2}) can be used, in order to get an upper-bound on ${\bf k}$, we require an upper-bound on the maximum absolute eigenvalue of $A$.  Due to the sparsity and diagonal dominance of $A$, a good estimate of the upper-bound can be found using the Gershgorin circle theorem (Def.~\ref{d:gersh} and Thm.~\ref{t:gersh}).
\begin{deff}[Ref.~\cite{MATRIX_BOOK}]
\label{d:gersh}
Let $A$ be a square complex matrix. Around every element $a_{ii}$ on the diagonal of the matrix, a circle with radius equal to the sum of the norms of the other elements in the same row ($\sum_{j\neq i}|a_{ij}|$) is known as a \emph{Gershgorin disc}.
\end{deff}
\begin{thm}[Ref.~\cite{MATRIX_BOOK}]
\label{t:gersh}
Every eigenvalue of a square complex matrix A lies in one of its Gershgorin discs.
\end{thm}
Since every eigenvalue must be contained in a Gershgorin disk, by finding the maximum absolute value of the limits of the disks will yield an upper-bound on the maximum modulus of the eigenvalues of $A$.

In the one-dimensional LSE case with no external potential and periodic boundary conditions, the $A$ matrix becomes circulant as defined by Def.~\ref{d:circ}.  In this case, the eigenvalues can be computed analytically by Thm.~\ref{t:circ}.  The upper-bound is then taken by finding the limit of the maximum eigenvalue as the size of the matrix goes to infinity.
\begin{deff}[Ref.~\cite{STB_circulant}]
\label{d:circ}
A circulant matrix is a square $N\times N$ matrix $C$ that can be fully specified by one vector, $\vec c = \{c_0,c_1,...,c_{N-1}\}$, which appears as the first column of $C$. The remaining columns of $C$ are each cyclic permutations of the vector with the offset equal to the column index.
\end{deff}
\begin{thm}[Ref.~\cite{STB_circulant}]
\label{t:circ}
The eigenvalues of a circulant matrix are given by
\[ \lambda_j = c_0 + c_{N-1}\,\omega_j + c_{N-2}\,\omega_j^2 + ... + c_1\,\omega_j^{N-1}, \qquad j = 0,..., {N-1},
\]
where 
\[
\omega_j = \exp \left(\frac{2\pi\, i\, j}{N}\right).
\]
\end{thm}

\section{Boundary conditions}
\label{s:bc}
Since boundary conditions of the spatial differencing in a PDE like the NLSE have the potential to alter the stability of a scheme, it is necessary to have stability results for each specific boundary condition one would like to use.  In this paper we limit ourselves to four boundary conditions which we feel are a good combination of simplicity and usefulness.  These boundary conditions are: periodic, Dirichlet, Laplacian-zero, and Modulus-Squared-Dirichlet.  As notation, we use the subscript $b$ to represent any boundary point, and $b-1$ to represent the grid position one point inward from the boundary in the normal direction.

For use with the stability analysis, it is desirable to formulate each boundary condition in terms of the temporal derivative in the form
\begin{equation}
\label{bb}
\left. \frac{\partial \Psi}{\partial t}\right|_b = \frac{i\,a}{h^2} B_b \Psi_b,
\end{equation}
and in terms of the spatial Laplacian in the form
\begin{equation}
\label{db}
\nabla^2\Psi_b = \frac{1}{h^2} D_b \Psi_b,
\end{equation}
where $B_b$ and $D_b$ are assumed to be real-valued constants (possibly differing per boundary point) and defined based on the specific boundary condition being used.  For periodic boundary conditions (or linear one-sided conditions not discussed here), these forms are not applicable.  Writing the boundary conditions in the forms of Eq.~(\ref{bb}) and Eq.~(\ref{db}) allows them to be expressed in the $A$ matrix as a single real-valued entry ($B_b$), and in the case of the form of the fourth-order differencing chosen here (see Sec.~\ref{s:1d2shocstb}), the near-boundary interior points will contain $D_b$ in their formulation.

\subsection{Periodic}
For periodic boundary conditions, any element of the scheme that is too small or too large in index (i.e. they are `off the grid') are simply replaced by the grid points on the opposite side of the grid.  In the case of the NLSE, periodic boundary conditions can be problematic especially in background-density situations due to the unpredictable phase jump from one side of the grid to the other.

\subsection{Dirichlet}
Dirichlet boundary conditions are defined as
\[
\Psi_b = C,
\]
where $C$ is a constant.  In terms of the temporal derivative of the NLSE, this condition is
\[
\left. \frac{\partial \Psi}{\partial t}\right|_b = 0,
\] 
in which case $B_b = 0$ in Eq.~(\ref{bb}).  When inserted into the NLSE, this condition in terms of the Laplacian is given by
\[
\nabla^2\Psi_b = -\frac{1}{a} (s|\Psi_b|^2 - V_b) \Psi_b,
\]
and therefore $D_b=-h^2/a (s|\Psi_b|^2 - V_b)$ in Eq.~(\ref{db}).

\subsection{Modulus-squared Dirichlet}
In some situations Dirichlet boundary condition can fail.  Such failure typically occurs in  simulations with a constant-density background, i.e. a constant value of $|\Psi|^2$ at the boundaries.  A standard Dirichlet condition will not work in such cases because it does not take into account the phase rotation of $\Psi$.  Instead, one would like to have the modulus-squared of the wavefunction to be constant at the boundaries, i.e. 
\[
|\Psi_b|^2 = C,
\]
where $C$ is a constant.  We have recently formulated a method for such a boundary condition (which is almost as easy to implement as Dirichlet) called the modulus-squared Dirichlet boundary condition \cite{ME_MSD}.  The MSD boundary condition is given in terms of the temporal derivative of the NLSE as
\begin{equation}
\label{msdfixedstb}
\Psi_{t,b} \approx i\,\mbox{Im}\left[\frac{\Psi_{t,b-1}}{\Psi_{b-1}}\right]\,\Psi_b.
\end{equation}
where $\partial \Psi_{b-1}/\partial t$ is computed by the interior scheme first, and then used to compute the boundary values.  Using the MSD boundary condition gives $B_b = (h^2/a)\mbox{Im}\left[\left.\frac{\partial \Psi}{\partial t}\right|_{b-1}\frac{1}{\Psi_{b-1}}\right]$, which is nonlinear, and not a constant independent of $\Psi$.  As shown in Ref.~\cite{ME_MSD}, due to the underlying assumptions of the MSD boundary condition, Eq.~(\ref{msdfixedstb}) can be viewed as  
\[
\left. \frac{\partial \Psi}{\partial t}\right|_b \approx i\,\Omega_{b-1} \Psi_b,
\]
where $\Omega_{b-1}$ is the real-valued frequency of the solution near the boundary.  Thus, $B_b$ would have the form $B_b = (h^2/a)\Omega_{b-1}$.  Therefore, we can linearize the MSD boundary condition by treating the $B_b$ term as a constant (which can change over the course of the simulation, similar to the nonlinearity of the NLSE).

When inserted into the NLSE, the MSD boundary condition of Eq.~(\ref{msdfixedstb}) yields
\begin{equation}
\label{msdlap}
\nabla^2\Psi_b \approx \left[ \mbox{Im}\left(i\, \frac{\nabla^2\Psi_{b-1}}{\Psi_{b-1}}\right) + \frac{1}{a}\,\left(N_{b-1} - N_b\right)\right]\Psi_b,
\end{equation}
where
\begin{equation}
N_b = s\,|\Psi_b|^2 - V_b, \qquad N_{b-1} = s\,|\Psi_{b-1}|^2 - V_{b-1}.
\end{equation}
and therefore $D_b = h^2\left[\mbox{Im}\left(i\, \frac{\nabla^2\Psi_{b-1}}{\Psi_{b-1}}\right) + \frac{1}{a}\,\left(N_{b-1} - N_b\right)\right]$.  This too is a nonlinear, non-constant term, and so must be treated as a constant in the same manner as the nonlinearity of the NLSE. 

\subsection{Laplacian-zero boundary condition}
The Laplacian-zero boundary condition is defined by 
\[
\nabla^2\Psi_b = 0,
\]
and therefore $D_b=0$.  In terms of the time-derivative of the NLSE, the L0 boundary condition is given as
\[
\left. \frac{\partial \Psi}{\partial t}\right|_b = i\,(s|\Psi_b|^2 - V_b)\,\Psi_b,
\]
making $B_b = (h^2/a)(s|\Psi_b|^2 - V_b)$.  This condition is as easy to implement as the Dirichlet, and can be useful in many situations.

To assist the stability analysis, a summary of the values of $B_b$ and $D_b$ for all the mentioned boundary conditions are given in Table.~\ref{t:bc} for future reference.  
\begin{table}[htbp] 
\centering 
\caption{Boundary condition terms for use with stability analysis.}
\begin{tabular}{|l|c|c|} \hline
Boundary Condition & $B_b$ & $D_b$ \\ \hline
& & \\
Dirichlet      & $0$                                      & $\dfrac{h^2}{a}(V_b-s|\Psi_b|^2)$ \\
\; &\; &\; \\
Laplacian-zero & $\dfrac{h^2}{a}(s|\Psi_b|^2-V_b)$         & $0$ \\
\; &\; &\; \\
MSD            & $\dfrac{h^2}{a}\mbox{Im}\left[\dfrac{\Psi_{t,b-1}}{\Psi_{b-1}}\right]$ & $h^2\left[\mbox{Im}\left(i\, \dfrac{\nabla^2\Psi_{b-1}}{\Psi_{b-1}}\right) + \dfrac{1}{a}\,\left(N_{b-1} - N_b\right)\right]$ \\ 
& & \\
\hline
\end{tabular}
\label{t:bc}
\end{table}
Many other boundary conditions exist for simulating the NLSE, in which case the analysis shown in this paper can be adapted to the other boundary conditions. 

\section{One-dimensional stability analysis}
\label{s:1dstb}
In the one-dimensional cases we analyze all four boundary conditions mentioned in Sec.~\ref{s:bc}.  As stated, periodic boundary conditions yield a matrix where (in the linear case with $s=0$ and $\vec V({\bf r})=0$) the eigenvalues can be computed analytically.  This allows the results obtained using the upper-bound methods (which we use with other boundary conditions) to be compared with the true eigenvalues giving an idea of how accurate they are.

\subsection{Second-order central difference}
\label{s:1dcdstb}
The second-order central difference is one dimension is given by
\[
\nabla^2\Psi_i = \left. \frac{\partial^2 \Psi}{\partial x^2} \right|_i \approx \frac{\Psi_{i+1} -2\Psi_i +\Psi_{i-1}}{h^2},
\]
and when implemented into the $A$ matrix, forms a matrix which is tridiagonal (except for the two boundary condition rows).

\subsubsection{Periodic boundary conditions}
\label{s:percd}
In order to obtain analytic expressions for the eigenvalues of $A$, we start with the LSE case with no external potential and periodic boundary conditions.  This yields the matrix
\[
A = \left[\begin {array}{cccccc} -2& 1& 0& 0& 1\\\noalign{\medskip}
                                  1&-2& 1& 0& 0\\\noalign{\medskip}
                                  0& \ddots&\ddots & \ddots & 0\\\noalign{\medskip}
                                  0& 0& 1 &-2& 1\\\noalign{\medskip}
                                  1& 0& 0& 1&-2\end {array} \right],
\]
which, as per Def.~\ref{d:circ}, is a circulant matrix with $\vec c = \{-2,1,0,...,0,1\}$.  Also, since $A$ is a real-valued symmetric matrix, by Thm.~\ref{t:real}, all eigenvalues are real and therefore the stability criteria of Eq.~(\ref{rk4stbeasy2}) can be used.  By Thm.~\ref{t:circ}, the eigenvalues of $A$ are given by
\[
\lambda_j = -2 + \exp\left[\frac{2\pi i j }{N}\right] +  \exp\left[\frac{2\pi i j (N-1)}{N}\right], \qquad j\in\{0,...N-1\}.
\]
The maximum value of $|\lambda_j|$ occurs either at $j=N/2$ if $N$ is even, or $j=(N\pm 1)/2$ if $N$ is odd. For $N$ even-valued we have
\[
|\lambda|_{\max} = \left|-2 + \exp\left[\pi i\right] + \exp\left[\pi i\right]^{N-1}\right|,
\]
which yields
\[
|\lambda|_{\max} = 4.
\]
For $N$ odd-valued we have
\[
|\lambda|_{\max} = \left|-2 - (-1)^{1/N} + (-1)^N\,(-1)^{-1/N}\right|,
\]
which yields
\[
|\lambda|_{\max} = \left|-2 - 2\,\cos\left(\frac{\pi}{N}\right)\right|.
\]
Taking $N\rightarrow \infty$, the maximum bound on the maximum absolute eigenvalue becomes
\[
|\lambda|_{\max} < 4.
\]
We therefore have an upper bound on the maximum absolute eigenvalue which, for even-valued $N$, is guaranteed to \emph{be} one of the eigenvalues.  The stability criteria of Eq.~(\ref{rk4stbeasy2}) is then formulated as
\begin{equation}
\label{1dcircbound}
{\bf k} < \frac{\sqrt{8}}{4}\frac{h^2}{a}.
\end{equation}

In the general case where $s \ne 0$ and/or $V({\bf r}) \ne 0$, the $A$ matrix is no longer circulant (since the values of the nonlinearity or external potential vary over the diagonal of $A$).  To get a bound on the maximum absolute eigenvalue, we make use of Thm.~\ref{t:gersh}.  The matrix $A$ has $N$ Gershgorin disks, since each diagonal entry of $A$ can be unique, but each disk has the same radius ($r=2$).  Also, since the diagonal entries can in theory take on any value, all Gershgorin disk limits must be examined.  This yields the stability bound
\begin{equation}
\label{1dcircboundL}
{\bf k} < \frac{\sqrt{8}}{\max\{\lVert\vec L\rVert_{\infty},\lVert \vec L-4\rVert_{\infty}\}}\,\frac{h^2}{a},
\end{equation}
where we have defined the elements of $\vec L$ to be 
\begin{equation}
\label{Leq}
L_i = \frac{h^2}{a}(s\,|\Psi_i|^2 - V_i),
\end{equation}
where the index $i$ spans over the entire grid.  It is important to note that all values of $\vec L$ are $O(h^2)$.   Thus, for $h\ll 1$, and reasonable values of $|\Psi|^2$ and $\vec V$, the linear bound of Eq.~(\ref{1dcircbound}) should be very close to the true bound of the nonlinear problem.

If we set $\vec L=0$ in Eq.~(\ref{1dcircboundL}), we recover the bound in Eq.~(\ref{1dcircbound}).  This shows that (in this case at least), using the Gershgorin circle theorem yields the true bound on the eigenvalues of $A$.

\subsubsection{Dirichlet, MSD, and L0 boundary conditions}
\label{s:1dcdbc2}
As shown in Sec.~\ref{s:bc}, Dirichlet, Laplacian-zero, and modulus-squared Dirichlet boundary conditions can all be viewed as single entries in the boundary value rows of the $A$ matrix, denoted as $B_b$.  As shown there, the values of $B_b$ are real-valued and their values for each boundary condition were given in Table~\ref{t:bc}.  Using such a formulation, the $A$ matrix becomes
\[
A = \left[ \begin {array}{cccccc}  B_0&0     &0  &0  &0         \\\noalign{\medskip}
                                   1  &L_1-2 &1  &0  &0         \\\noalign{\medskip}
                                   0  &\ddots   &\ddots &\ddots          &0\\\noalign{\medskip}
                                   0  &0      &1  &L_{N-2}-2 &1\\\noalign{\medskip}
                                   0  &0       &0  &0         &B_{N-1}\end {array} \right]. 
\]
In order to use the simple stability criteria of Eq.~(\ref{rk4stbeasy2}), we once again need to show that all eigenvalues of $A$ are purely real.  The $A$ matrix is no longer symmetric, however it is easy to see that $B_0$ and $B_{N-1}$ are eigenvalues of $A$, and the remaining eigenvalues of $A$ are equivalent to the eigenvalues of the matrix $A^{\prime}$ defined as
\[
A^\prime = \left[\begin {array}{cccccc}
L_1-2&1    &0  &0        &0\\\noalign{\medskip}
1    &L_2-2&1  &0        &0\\\noalign{\medskip}
0    &\ddots    &\ddots &\ddots       &0\\\noalign{\medskip}
0    &0    &1 &L_{N-3}-2&1\\\noalign{\medskip}
0    &0    &0  &1        &L_{N-2}-2
\end {array} \right].
\]
Since $A^{\prime}$ is real-valued symmetric, we can use the stability bound of Eq.~(\ref{rk4stbeasy2}).  

We now need to find an upper bound on the absolute value of the eigenvalues of $A^{\prime}$.  We use the Gershgorin circle theorem to find all unique Gershgorin disks and take the limits of the disks to find the bounds on the absolute eigenvalues. Many of the Gershgorin disks are similar, differing only in the value of $L_i$ of the specific row.  Therefore, each disk of different centers and radii has a subset of $L_i$ values relevant to it.  Although in the current one-dimensional setting it is simple to define the subsets, in higher-dimensional settings, it can become burdensome to separate out each subset of $\vec L$ relevant to each Gershgorin disk of the same center and radius.  Therefore, for practicality purposes, we define our bounds using all possible values of $L_i$ for each Gershgorin disk center and radius.  This may make the resulting stability bound slightly higher than necessary in certain cases, but this is outweighed by the ease-of-use of the simplified bounds.  The unique forms of the Gershgorin disks of $A^{\prime}$ are shown in Table~\ref{t:1dcdgd}.
\begin{table}[htbp] 
\centering 
\caption{Unique forms of the Gershgorin disk centers ($a_{ii}$) and radii ($r_i$) for the $A^{\prime}$ matrix of the one-dimensional second-order central difference scheme.}
\begin{tabular}{|c|c|} \hline
$a_{ii}$    & $r_i=\sum_{i\ne j} |a_{ij}|$ \\ \hline
$L_i - 2$   & $1$  \\
$L_i - 2$   & $2$  \\
\hline
\end{tabular}
\label{t:1dcdgd}
\end{table}
The resulting general stability bounds are
\begin{equation}
\label{kbg}
{\bf k} < \frac{\sqrt{8}}{\max\{\lVert\vec B\rVert_{\infty},\lVert \forall L_i, L_i - \vec G\rVert_{\infty}\}}\,\frac{h^2}{a},
\end{equation}
where $\vec B$ are all boundary condition values (in this case $B_0$ and $B_{N-1}$), and $\vec G$ is defined as
\begin{equation}
\label{G1DCD}
\vec G = \left\{4,3,1,0 \right\}.
\end{equation}
In general, all possible values of $\vec  G$ must be taken into consideration since there is no theoretical restriction on what values $\vec L$ can take.  However, in certain specific circumstances, some of the values of $\vec G$ can be ignored (for example, when $s \le 0$ and $V({\bf r}) \ge 0$, only the largest magnitude value in $\vec G$ is needed).

\subsection{Fourth order central difference}
\label{s:1d2shocstb}
The standard fourth order central difference scheme is given by
\begin{equation}
\label{1d4cd}
\nabla^2\Psi_i = \left. \frac{\partial^2 \Psi}{\partial x^2} \right|_i \approx \frac{-\Psi_{i+2} +16\Psi_{i+1} -30\Psi_i +16\Psi_{i-1} - \Psi_{i-2}}{12\,h^2}.
\end{equation}
The stability analysis follows directly from the second-order case.  The only major difference is that since the fourth order stencil is five points wide, the grid points near the boundary may need special consideration for the different boundary conditions.  For our purposes here, we use the two-step high-order compact (2SHOC) version of the fourth-order scheme as described in Ref.~\cite{ME_2SHOC}, in which case the near-boundary points can be formulated by combining the two steps of the 2SHOC scheme.

\subsubsection{Periodic boundary condition}
\label{s:per2shoc}
In the periodic case, no special attention is needed near the boundaries, and the $A$ matrix in the LSE case with no external potential is
\[
A = \left[\begin {array}{cccccccc}
-15/6 &   4/3 & -1/12 &     0 &     0 &     -1/12 &   4/3\\\noalign{\medskip}
  4/3 & -15/6 &   4/3 & -1/12 &     0 &         0 & -1/12\\\noalign{\medskip}
-1/12 &   4/3 & -15/6 &   4/3 & -1/12 &         0 &     0\\\noalign{\medskip}
    0 &   \ddots &   \ddots &   \ddots &  \ddots &  \ddots &         0\\\noalign{\medskip}
    0 &         0 & -1/12 &   4/3 & -15/6 &   4/3 & -1/12\\\noalign{\medskip}
-1/12 &     0 &        0 & -1/12 &   4/3 & -15/6 &   4/3\\\noalign{\medskip}
  4/3 & -1/12 &        0 &     0 & -1/12 &   4/3 & -15/6\end {array} \right],
\]
which is a circulant matrix, and its eigenvalues are therefore
\begin{alignat}{2}
\lambda_j = &-\frac{15}{6} + \frac{4}{3}\,\exp\left[\frac{2\pi i j }{N}\right] - \frac{1}{12}\,\exp\left[\frac{4\pi i j }{N}\right]  \notag \\
&- \frac{1}{12}\,\exp\left[\frac{2(N-2)\pi i j }{N}\right] + \frac{4}{3}\,\exp\left[\frac{2(N-1)\pi i j }{N}\right]. \notag
\end{alignat}
The maximum absolute value once again occurs occurs at either $j=N/2$ if $N$ is even, or $j=(N\pm 1)/2$ if $N$ is odd. For $N$ even-valued we have
\[
\lambda_{N/2} = -\frac{15}{6} - \frac{4}{3} - \frac{1}{12} - \frac{1}{12}(-1)^{N-2} + \frac{4}{3}(-1)^{N-1} = -\frac{16}{3}.
\]
For $N$ odd, we have
\[
\lambda_{(N+1)/2} = -\frac{15}{6} - \frac{4}{3}\,\left((-1)^{1/N} + (-1)^{-1/N}\right) - \frac{1}{12}\left((-1)^{2/N} + (-1)^{-2/N}\right),
\]
which yields
\[
\lambda_{(N+1)/2} = -\frac{15}{6} -\frac{4}{3}\,\left(2\,\cos\left(\frac{\pi}{N}\right)\right) - \frac{1}{12}\,\left(2\,\cos\left(\frac{2\pi}{N}\right)\right).
\]
As $N\rightarrow \infty$, $|\lambda| \rightarrow \frac{16}{3}$, which is the same bound as the $N$-even case.  Thus, the stability bound is given by
\begin{equation}
\label{klinhoc1d}
{\bf k} < \left(\frac{3}{4}\right) \frac{\sqrt{8}}{4} \frac{h^2}{a},
\end{equation}
which we note is only a $25\%$ reduction of the second-order bound of Eq.~(\ref{1dcircbound}).  In the general case where $\vec L \ne 0$, the bound becomes
\begin{equation}
{\bf k} < \frac{6\sqrt{2}}{\max\{\lVert3\vec L-16\rVert_{\infty},\lVert 3\vec L+1\rVert_{\infty}\}} \frac{h^2}{a}.
\end{equation}
If $s \le 0$ and $V({\bf r}) \ge 0$, the first term in the denominator is the maximum of the two terms, and the resulting stability bound is equivalent to that found in Ref.~\cite{RK4_2CNLSE_STB} for using the RK4 scheme and fourth-order spatial differencing with the coupled NLSE.

\subsubsection{Dirichlet, MSD, and L0 boundary conditions}
\label{s:1d2shoc}
As per Sec.~\ref{s:bc}, we formulate all three boundary conditions in terms of a $B_b$ entry in the $A$ matrix.  As discussed, an important issue is that we need to handle the grid points near the boundary due to the width of the scheme.  A common way of dealing with the closest-interior points is to compute the Laplacian at those points using second-order differencing, however this can lead to the overall scheme becoming second-order.  However, since we are using the 2SHOC version of the fourth-order differencing, we can derive the closest-interior points which, if the assumptions of the chosen boundary conditions hold, should maintain fourth-order accuracy.  In one dimension, the 2SHOC scheme is defined as  \cite{ME_2SHOC}
\begin{alignat}{3}
 &1) \qquad &D_i &= \frac{1}{h^2}\left(\Psi_{i+1} - 2\Psi_i + \Psi_{i-1}\right), \label{2shoc1dstb} \\
 &2) \qquad &\nabla^2\Psi_i &\approx \frac{7}{6}D_i - \frac{1}{12}\left(D_{i+1} + D_{i-1}\right). \label{2shoc1d2stb}
\end{alignat}
In the first step, the second-order Laplacian is computed with the chosen boundary condition applied to it.  Next, the result is used to compute the fourth-order Laplacian. As mentioned in Ref.~\cite{ME_2SHOC}, this two-step scheme is equivalent to the standard wide-stencil of Eq.~(\ref{1d4cd}) for the interior points.  We use the form of Eq.~(\ref{db}) for the boundary conditions on the Laplacian, and after combining the steps of Eq.~(\ref{2shoc1dstb}) and Eq.~(\ref{2shoc1d2stb}), we get the matrix
\[
A = \left[\begin {array}{cccccccc}
B_0         & 0         & 0        & 0      & 0            & 0             & 0 \\\noalign{\medskip}
\frac{14-D_0}{12} & L_1-\frac{29}{12} & \frac{4}{3}     & -\frac{1}{12}  & 0            & 0             & 0 \\\noalign{\medskip}
-\frac{1}{12}       & \frac{4}{3}       & L_2-\frac{15}{6} & \frac{4}{3}    & -\frac{1}{12}        & 0             & 0 \\\noalign{\medskip}
0           & \ddots    & \ddots   & \ddots & \ddots       & \ddots        & 0 \\\noalign{\medskip}
0           & 0         &-\frac{1}{12}    & \frac{4}{3}    & L_{N-3}-\frac{15}{6} & \frac{4}{3}           & -\frac{1}{12} \\\noalign{\medskip}
0           & 0         & 0        & -\frac{1}{12}  & \frac{4}{3}          & L_{N-2}-\frac{29}{12} & \frac{14-D_{N-1}}{12} \\\noalign{\medskip}
0           & 0         & 0        & 0      & 0            & 0             & B_{N-1}
\end {array} \right],
\]
where the $B_b$ and $D_b$ ($b\in\{0,N-1\}$) terms for each different boundary condition are once again given in Table~\ref{t:bc}.  As in Sec.~\ref{s:1dcdbc2}, the $A$ matrix is not symmetric and has eigenvalues equal to $\vec B$.  The remaining eigenvalues are those of the matrix $A^{\prime}$ defined as
\[
A^{\prime} = \left[\begin {array}{ccccccc}  
L_1-\frac{29}{12} & \frac{4}{3}      & -\frac{1}{12}    & 0      & 0            & 0            & 0 \\\noalign{\medskip}
\frac{4}{3}      & L_2-\frac{15}{6} & \frac{4}{3}       & -\frac{1}{12}  & 0            & 0            & 0 \\\noalign{\medskip}
-\frac{1}{12}     & \frac{4}{3}      & L_3-\frac{15}{6} & \frac{4}{3}    & -\frac{1}{12}        & 0            & 0 \\\noalign{\medskip}
0         & \ddots   & \ddots   & \ddots & \ddots       & \ddots       & 0 \\\noalign{\medskip}
0         & 0        & -\frac{1}{12}    & \frac{4}{3}     & L_{N-4}-\frac{15}{6} & \frac{4}{3}          & -\frac{1}{12} \\\noalign{\medskip}
0         & 0        & 0        & -\frac{1}{12}  & \frac{4}{3}           & L_{N-3}-\frac{15}{6} & \frac{4}{3} \\\noalign{\medskip}
0         & 0        & 0        & 0      & -\frac{1}{12}        & \frac{4}{3}          & L_{N-2}-\frac{29}{12}
\end {array} \right],
\]
which is once again real-symmetric so the bounds of Eq.~(\ref{rk4stbeasy2}) can be used.  It is interesting to note that the values of $D_b$ do not appear in any of the eigenvalues of $A^{\prime}$.

\begin{table}[htbp] 
\centering
\caption{Unique forms of the Gershgorin disk centers ($a_{ii}$) and radii ($r_i$) for the $A^{\prime}$ matrix of the one-dimensional fourth-order 2SHOC scheme.}
\begin{tabular}{|l|c|} \hline
$a_{ii}$    & $r_i=\sum_{i\ne j} |a_{ij}|$ \\ \hline
$L_i - 5/2$   & $11/4$  \\
$L_i - 5/2$   & $17/6$  \\
$L_i - 29/12$ & $17/12$ \\
\hline
\end{tabular}
\label{t:1d2shocgd}
\end{table}
The unique forms (see the discussion in Sec.~\ref{s:1dcdbc2}) of the Gershgorin disk centers and radii of $A^{\prime}$ are shown in Table~\ref{t:1d2shocgd}.  The full stability bound is the same form as Eq.~(\ref{kbg}), but with $\vec G$ defined as
\begin{equation}
\label{G1D2SHOC}
\vec G = \frac{1}{12} \times \left \{64,63,46,12,-3,-4\right\}.
\end{equation}
Once again, in the most general case, all values of Eq.~(\ref{G1D2SHOC}) must be considered in finding the maximum allowed time-step value.

\section{Two-dimensional stability analysis}
\label{s:2d}
In higher dimensions, the $A$ matrix is formed by unwrapping the solution into a one-dimensional vector and then formulating the scheme matrix accordingly.  

In Secs. \ref{s:percd} and \ref{s:per2shoc} we noted that the stability bounds given using the Gershgorin circle theorem were equivalent to those obtained analytically for the linear case with periodic boundary conditions.  We therefore justify relying exclusively on the Gershgorin theorem for higher dimensions, and focus on the stability bounds for Dirichlet, MSD, and Laplacian-zero boundary conditions (since the periodic boundary condition bounds will be a subset of the bounds computed for the other boundary conditions).  

\subsection{Second-order central difference}
\label{s:2dcdstb}
\begin{figure}[htbp]
\centering
\includegraphics[width=3.2in]{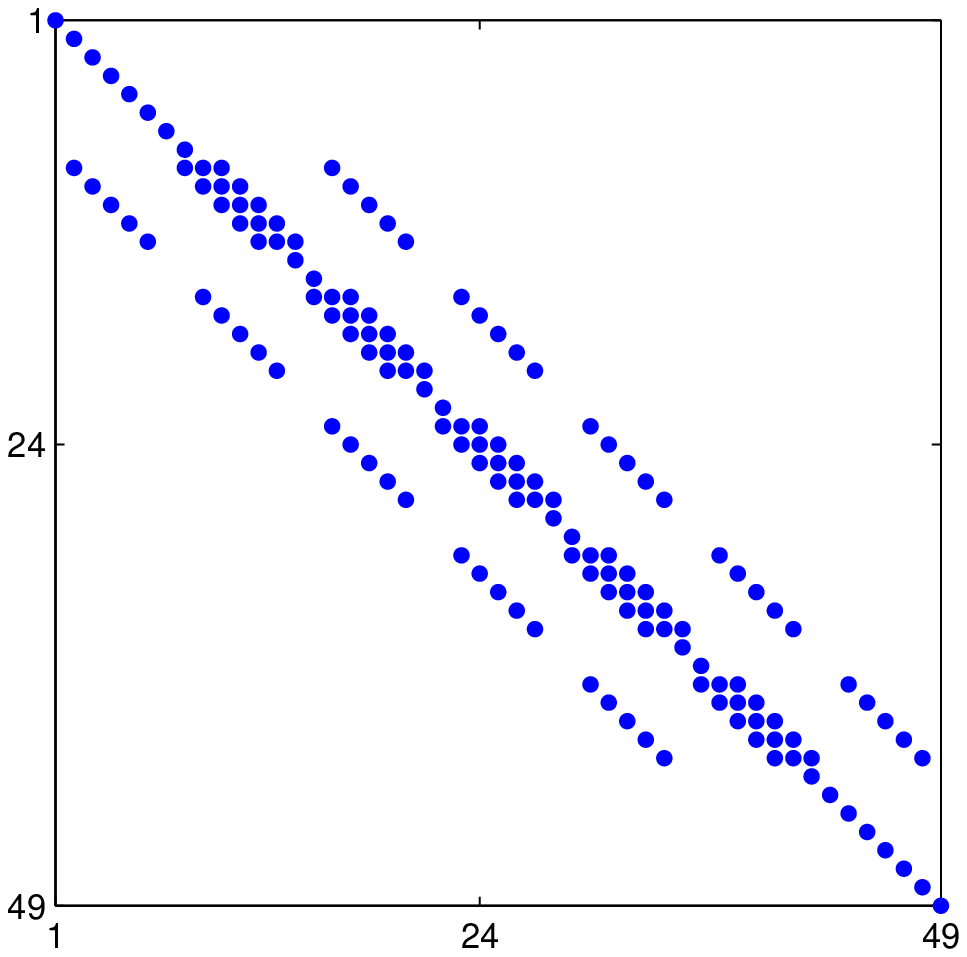}
\includegraphics[width=3.2in]{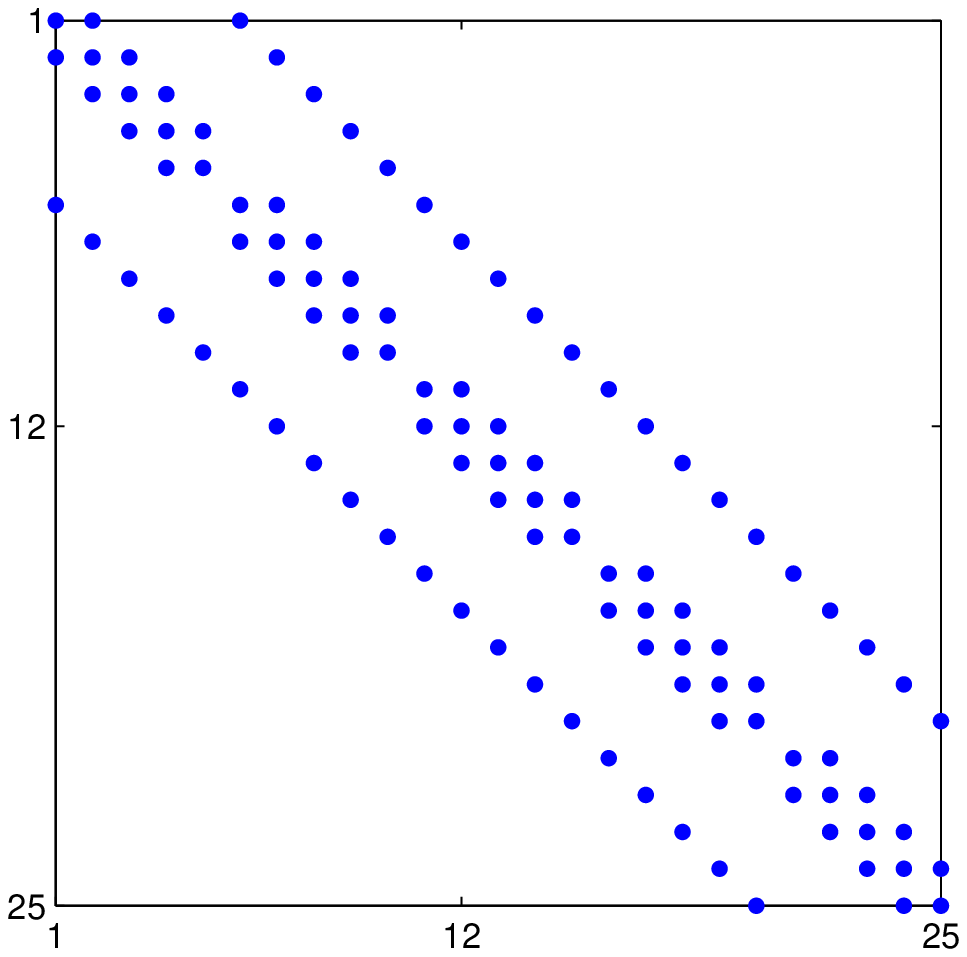}
\caption[Matrix structure for two-dimensional CD scheme.]{(Color online) Form of scheme matrix $A$ and $A^{\prime}$ for second-order central differencing of the two-dimensional NLSE.  The dots represent non-zero entries of the matrices.  The matrices shown are for a $7\times 7$ grid. \label{f:A2D}}
\end{figure}

The second-order central difference scheme in two dimensions is given by
\begin{equation}
\begin{tabular}{c} 
$\qquad \nabla^2 \Psi_{i,j} = \left. \dfrac{\partial^2 \Psi}{\partial x^2} \right|_{i,j} + \left. \dfrac{\partial^2 \Psi}{\partial y^2} \right|_{i,j} \approx \dfrac{1}{h^2}$
\begin{tabular}{|c|c|c|} \hline
  &  $1$ &   \\ \hline
$1$ & $-4$ & $1$ \\ \hline
  &  $1$ &  \\ \hline
\end{tabular}
$\Psi_{i,j}.$
\end{tabular}
\end{equation}
The corresponding $A$ matrix has a tri-banded structure, with diagonal sub-sections corresponding to the boundary values.  The form of the $A$ matrix is shown in Fig.~\ref{f:A2D} (we do not show the values of the entries of the matrix due to space considerations, but they can be obtained through symbolic math codes).  As in the one-dimensional case, all diagonal entries (which are the boundary value entries $B_b$) of $A$ are eigenvalues, and the remaining eigenvalues are real and equivalent to those of a matrix $A^{\prime}$ which is real-symmetric, thus allowing the use of the bounds in Eq.~(\ref{rk4stbeasy2}).  The form of $A^\prime$ is also shown in Fig.~\ref{f:A2D}.

The unique forms of the Gershgorin disk centers and radii for $A^{\prime}$ are shown in Table~\ref{t:2dcdgd}. 
\begin{table}[htbp] 
\centering 
\caption{Gershgorin disk centers and radii for the $A^{\prime}$ matrix of the two-dimensional central-difference scheme.}
\begin{tabular}{|c|c|} \hline
$a_{ii}$  & $r_i = \sum_{i\ne j} |a_{ij}|$ \\ \hline
$L_i - 4$ & $2$ \\
$L_i - 4$ & $3$ \\
$L_i - 4$ & $4$ \\
\hline
\end{tabular}
\label{t:2dcdgd}
\end{table}
The stability bounds are once again the same as in Eq.~(\ref{kbg}) with $\vec B$ being the set of all boundary values $B_b$ and with $\vec G$ now defined as
\begin{equation}
\label{G2DCD}
\vec G = \{0,1,2,6,7,8\}.
\end{equation}
In the linear case ($s=0$) with no external potential and periodic, Dirichlet or Laplacian-zero boundary conditions, we get the linear stability bound
\begin{equation}
\label{kbound2dcdlin}
{\bf k} < \frac{\sqrt{8}}{8}\,\frac{h^2}{a}.
\end{equation}
As before, since within the $A$ matrix, all boundary, potential, and nonlinear terms are $O(h^2)$ the simple bound of Eq.~(\ref{kbound2dcdlin}) with slight adjustment can be used in many applications.

\subsection{Fourth-order central difference}
\label{s:2d2shocstb}
\begin{figure}[htbp]
\centering
\includegraphics[width=3.2in]{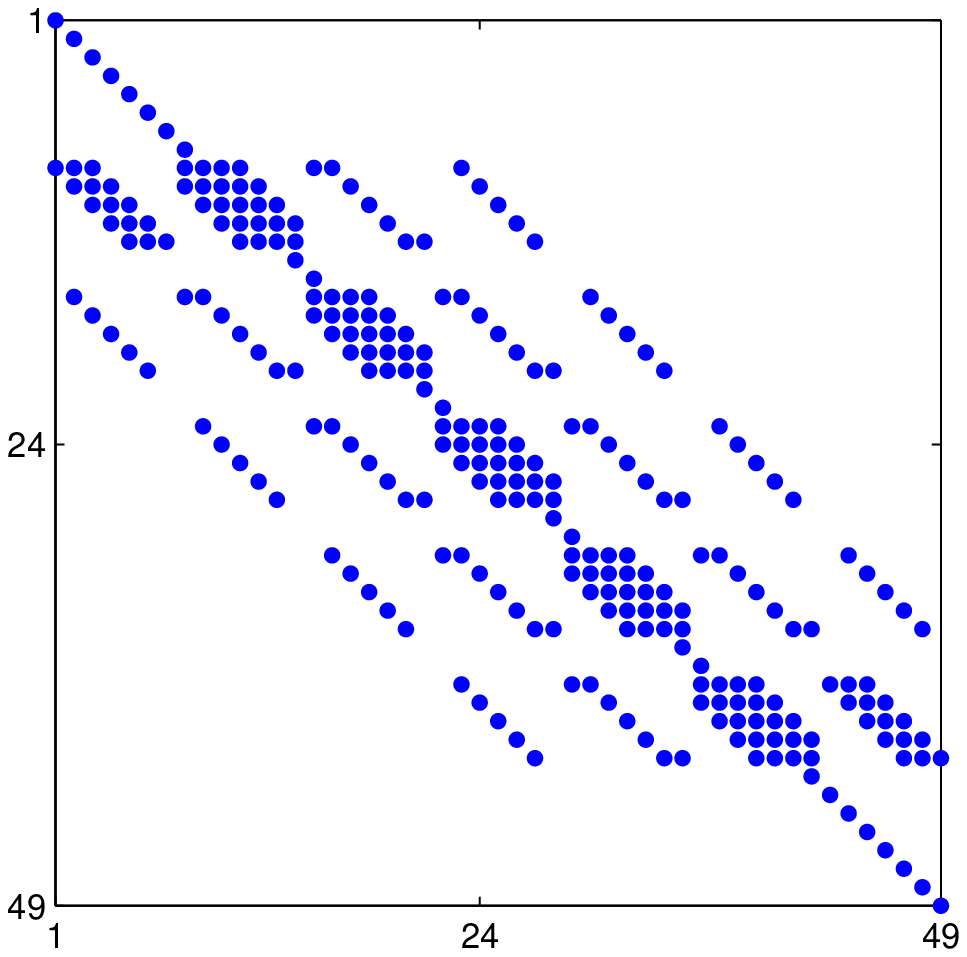}
\includegraphics[width=3.2in]{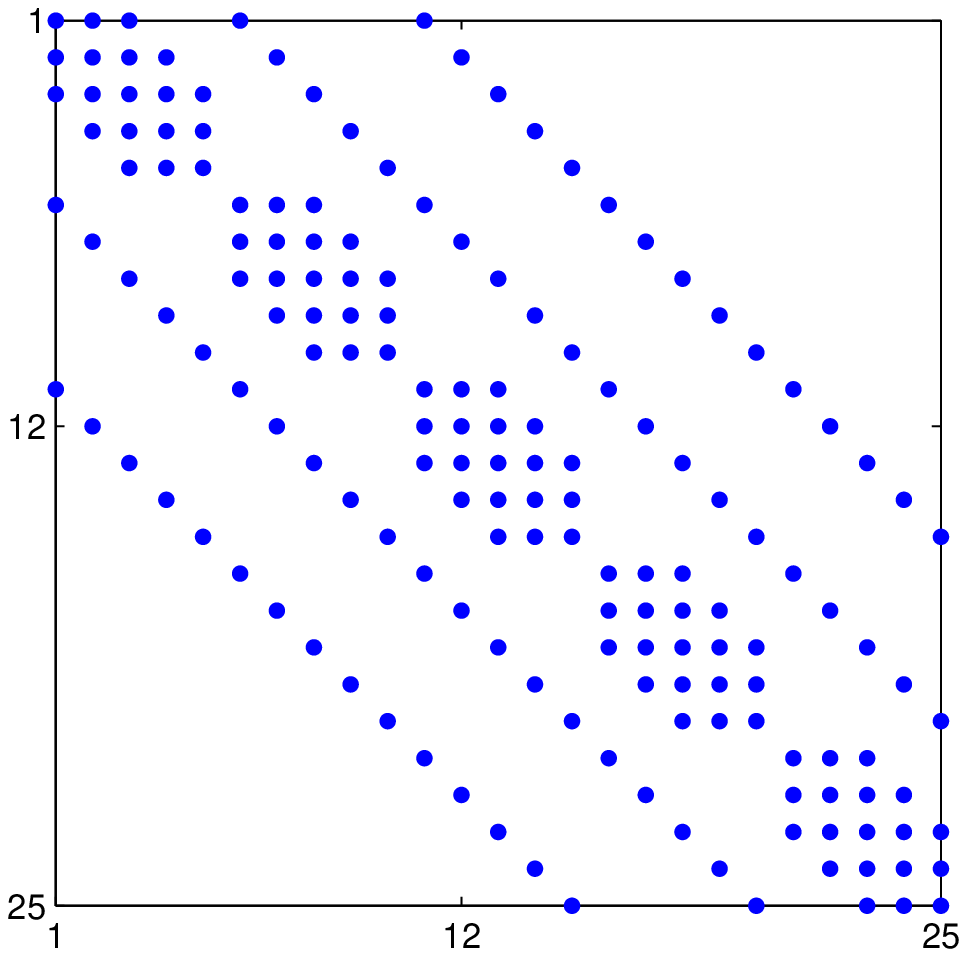}
\caption[Matrix structure for two-dimensional 2SHOC scheme.]{(Color online) Form of scheme matrix $A$ and $A^{\prime}$ for the fourth-order 2SHOC scheme of the two-dimensional NLSE.  The dots represent non-zero entries of the matrices.  The matrices shown are for a $7\times 7$ grid. \label{f:A2D2SHOC}}
\end{figure}

The fourth-order central difference scheme in two dimensions is given by 
\begin{equation}
\label{2dnonhoc}
\begin{tabular}{c} 
$\nabla^2\Psi_{i,j} \approx -\dfrac{1}{12\,h^2}$
\begin{tabular}{|c|c|c|c|c|} \hline
  &     &   $1$ &     &   \\ \hline
  &     & $-16$ &     &   \\ \hline
$1$ & $-16$ &  $60$ & $-16$ & $1$ \\ \hline
  &     & $-16$ &     &   \\ \hline
  &     &   $1$ &     &   \\ \hline
\end{tabular}
$\Psi_{i,j}.$
\end{tabular}
\end{equation}
The low-storage version of the 2SHOC equivalent scheme is defined as \cite{ME_2SHOC}
\begin{alignat}{3}
&\begin{tabular}{ll} 
$1)$ & $D_{i,j} = \dfrac{1}{h^2}$ 
\begin{tabular}{|c|c|c|} \hline
  &  $1$ &   \\ \hline
$1$ & $-4$ & $1$ \\ \hline
  &  $1$ &  \\ \hline
\end{tabular}
$\Psi_{i,j}$
\end{tabular}
\\
\;&\; \notag
\\
&\begin{tabular}{ll} 
$2)$ & $\nabla^2\Psi_{i,j} \approx -\dfrac{1}{12}$
\begin{tabular}{|c|c|c|} \hline
  &   $1$ &   \\ \hline
$1$ & $-12$ & $1$ \\ \hline
  &   $1$ &   \\ \hline
\end{tabular}
$D_{i,j} + \dfrac{1}{6\,h^2}$  
\begin{tabular}{|c|c|c|} \hline
$1$ &    & $1$ \\ \hline
  & $-4$ &   \\ \hline
$1$ &    & $1$ \\ \hline
\end{tabular}
$\Psi_{i,j}.$
\end{tabular}
\end{alignat}
The corresponding $A$ matrix has a five-banded structure.  The structure of the $A$ matrix and its corresponding $A^{\prime}$ matrix are shown in Fig.~\ref{f:A2D2SHOC}.

The resulting unique forms of the Gershgorin disk centers and radii for $A^{\prime}$ are shown in Table~\ref{t:2d2shocgd}. 
\begin{table}[htbp] 
\centering 
\caption{Gershgorin disk centers ($a_{ii}$) and radii ($r_i$) for the $A^{\prime}$ matrix of the two-dimensional 2SHOC scheme.}
\begin{tabular}{|l|c|} \hline
$a_{ii}$      & $r_i = \sum_{i\ne j} |a_{ij}|$ \\ \hline
$L_i - 5$     & $11/2$ \\
$L_i - 5$     & $67/12$\\
$L_i - 5$     & $17/3$ \\
$L_i - 29/6$  & $17/6$ \\
$L_i - 59/12$ & $25/6$ \\
$L_i - 59/12$ & $17/4$ \\
\hline
\end{tabular}
\label{t:2d2shocgd}
\end{table}
The linearized stability bound is once again Eq.~(\ref{kbg}) but with $\vec G$ defined as
\begin{equation}
\label{G2D2SHOC}
\vec G = \dfrac{1}{12} \times \left\{128,127,126,110,109,92,24,9,8,-6,-7,-8\right \}.
\end{equation}
The linear bound (with $s=0$ and $V({\bf r})=0$) is then given by
\begin{equation}
\label{kbound2d2shoclin}
{\bf k} < \frac{3\sqrt{8}}{32}\,\frac{h^2}{a},
\end{equation}
which, as in the one-dimensional case, is $25\%$ lower than the second-order linear bound given in Eq.~(\ref{kbound2dcdlin}).

\section{Three-dimensional stability analysis}
\label{s:3dstb}
For the stability analysis in three dimensions, the same procedure utilized in the two-dimensional case of Sec.~\ref{s:2d} is used.

\subsection{Second-order central difference}
\label{3dcdstb}
\begin{figure}[htbp]
\centering
\includegraphics[width=3.2in]{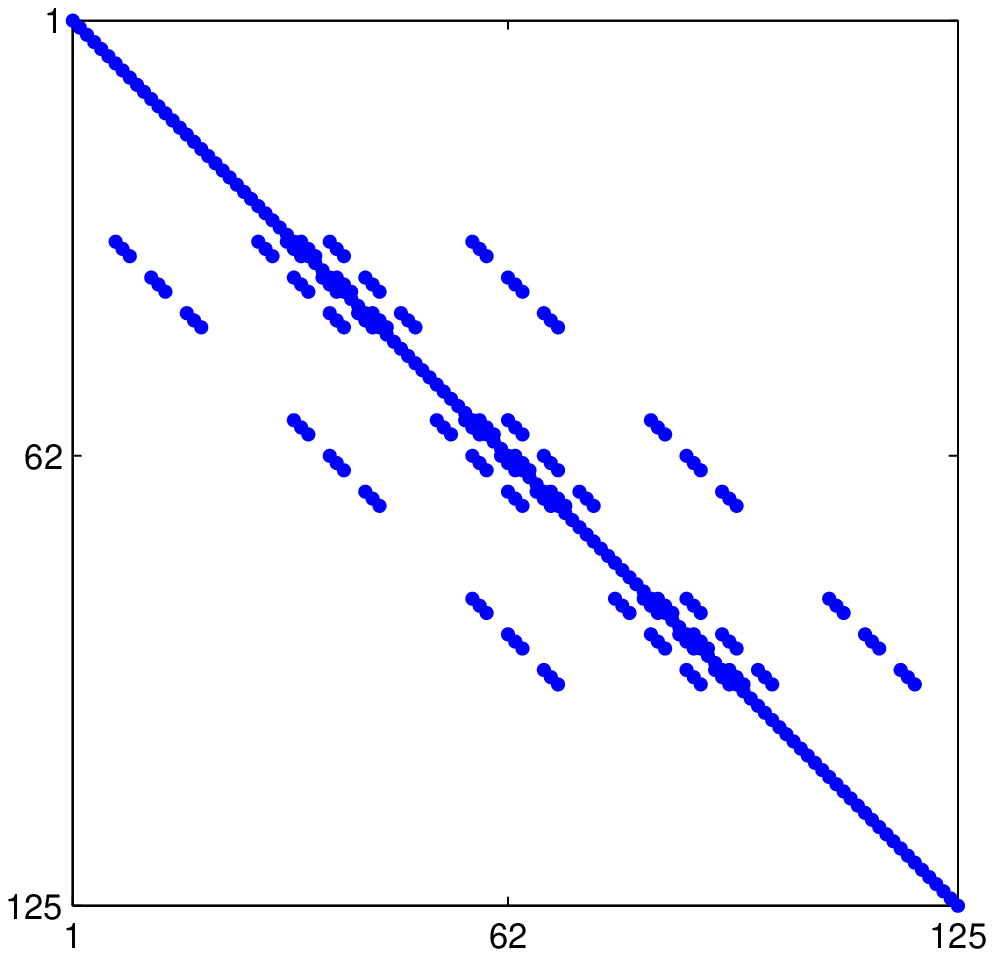}
\includegraphics[width=3.2in]{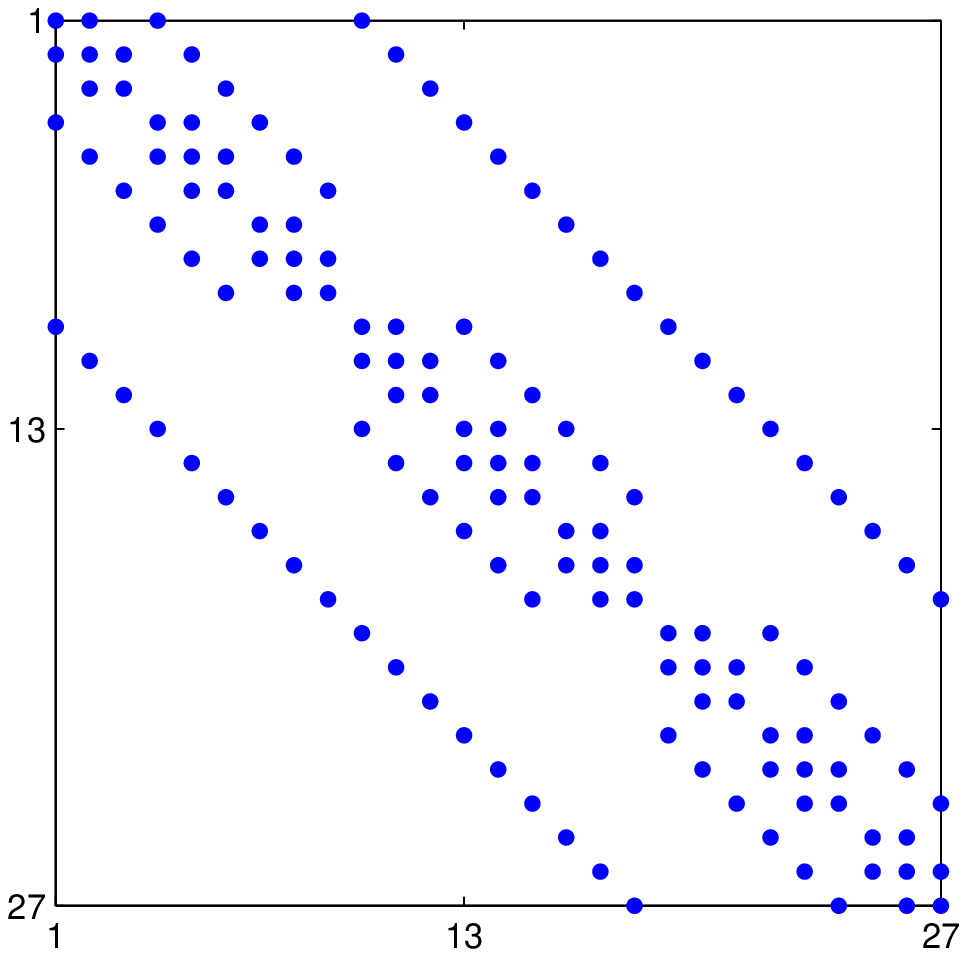}
\caption[Matrix structure for three-dimensional CD scheme.]{(Color online) Form of scheme matrix $A$ and $A^{\prime}$ for the second-order central difference scheme of the three-dimensional NLSE.  The dots represent non-zero entries of the matrices.  The matrices shown are for a $5 \times 5$ grid. \label{f:A3DCD}}
\end{figure}

The second-order central difference scheme in three dimensions is given by
\begin{equation}
\label{3dcds}
\begin{tabular}{ll}
$\nabla^2\Psi_{i,j,k}$ & $ = \left. \dfrac{\partial^2 \Psi}{\partial x^2} \right|_{i,j,k} + \left. \dfrac{\partial^2 \Psi}{\partial y^2} \right|_{i,j,k} + \left. \dfrac{\partial^2 \Psi}{\partial z^2} \right|_{i,j,k} $ \\
\; \; \\
\; & $\approx \dfrac{1}{h^2}$
$\left(\;
\begin{tabular}{|c|c|c|} \hline
  &   &   \\ \hline
  & $1$ &  \\ \hline
  &   &   \\ \hline
\end{tabular}
 \, \Psi_{i,j+1,k} +  
 \begin{tabular}{|c|c|c|} \hline
  & $1$ &   \\ \hline
$1$  & $-6$ & $1$ \\ \hline
  & $1$ &   \\ \hline
\end{tabular}
\, \Psi_{i,j,k} +  
\begin{tabular}{|c|c|c|} \hline
  &   &   \\ \hline
  & $1$ & \\ \hline
  &   &   \\ \hline
\end{tabular}
 \, \Psi_{i,j-1,k} \; \right),$ 
 \end{tabular}
\end{equation}
and the structure of the corresponding $A$ and $A^{\prime}$ matrix are given in Fig.~\ref{f:A3DCD}.
The unique forms of the Gershgorin disk centers ($a_{ii}$) and radii ($r_i$) for $A^{\prime}$ are shown in Table~\ref{t:3dcdgd}.
\begin{table}[htbp] 
\centering 
\caption{Gershgorin disk centers ($a_{ii}$) and radii ($r_i$) for the $A^{\prime}$ matrix of the three-dimensional central difference scheme.}
\begin{tabular}{|c|c|} \hline
$a_{ii}$   & $r_i = \sum_{i\ne j} |a_{ij}|$ \\ \hline
$L_i - 6$  & $3$ \\
$L_i - 6$  & $4$ \\
$L_i - 6$  & $5$ \\
$L_i - 6$  & $6$ \\
\hline
\end{tabular}
\label{t:3dcdgd}
\end{table}
The stability bounds of Eq.~(\ref{kbg}) in this case has $\vec G$ defined as
\begin{equation}
\label{G3DCD}
\vec G = \{12,11,10,9,3,2,1,0\}.
\end{equation}
In the linear case ($s=0$) with no external potential and periodic, Dirichlet or Laplacian-zero boundary conditions, the linear stability bound becomes
\begin{equation}
\label{kbound3dcdlin}
{\bf k} < \frac{\sqrt{8}}{12}\,\frac{h^2}{a}.
\end{equation}

\subsection{Fourth-order central difference}
\label{s:3d2shocstb}
\begin{figure}[htbp]
\centering
\includegraphics[width=3.2in]{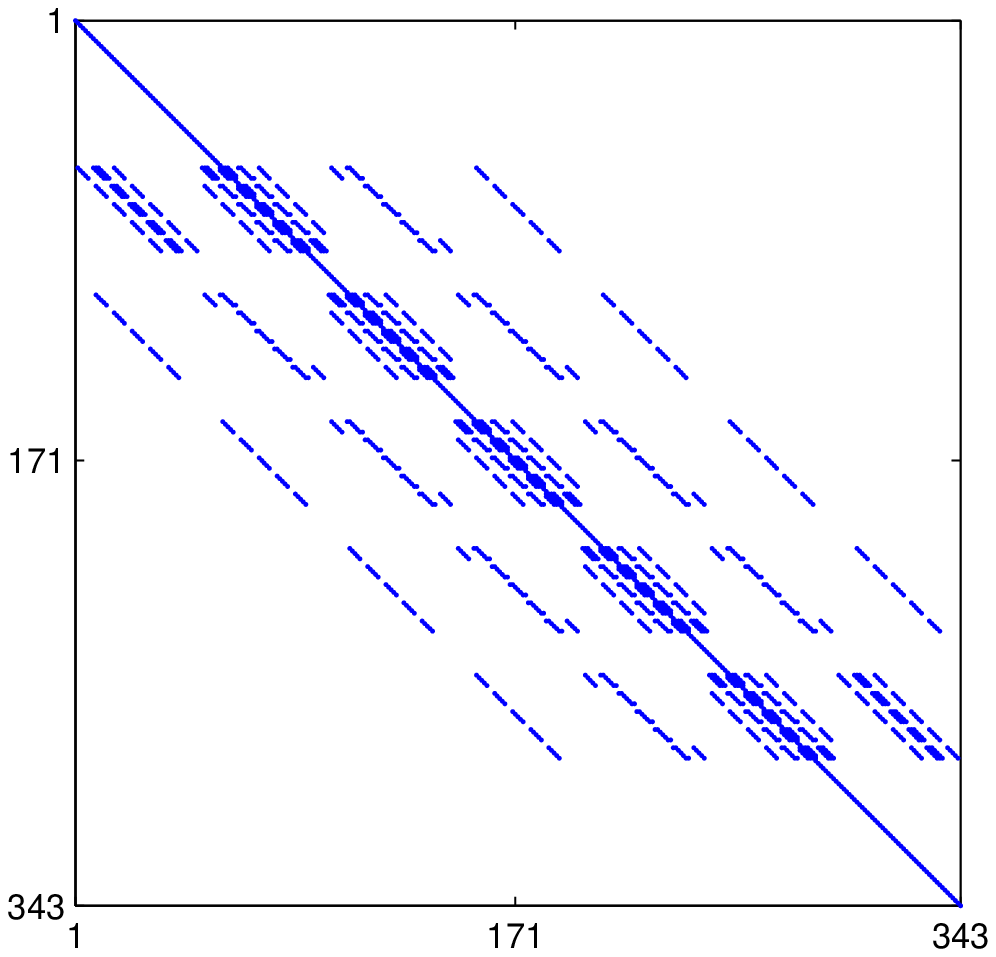}
\includegraphics[width=3.2in]{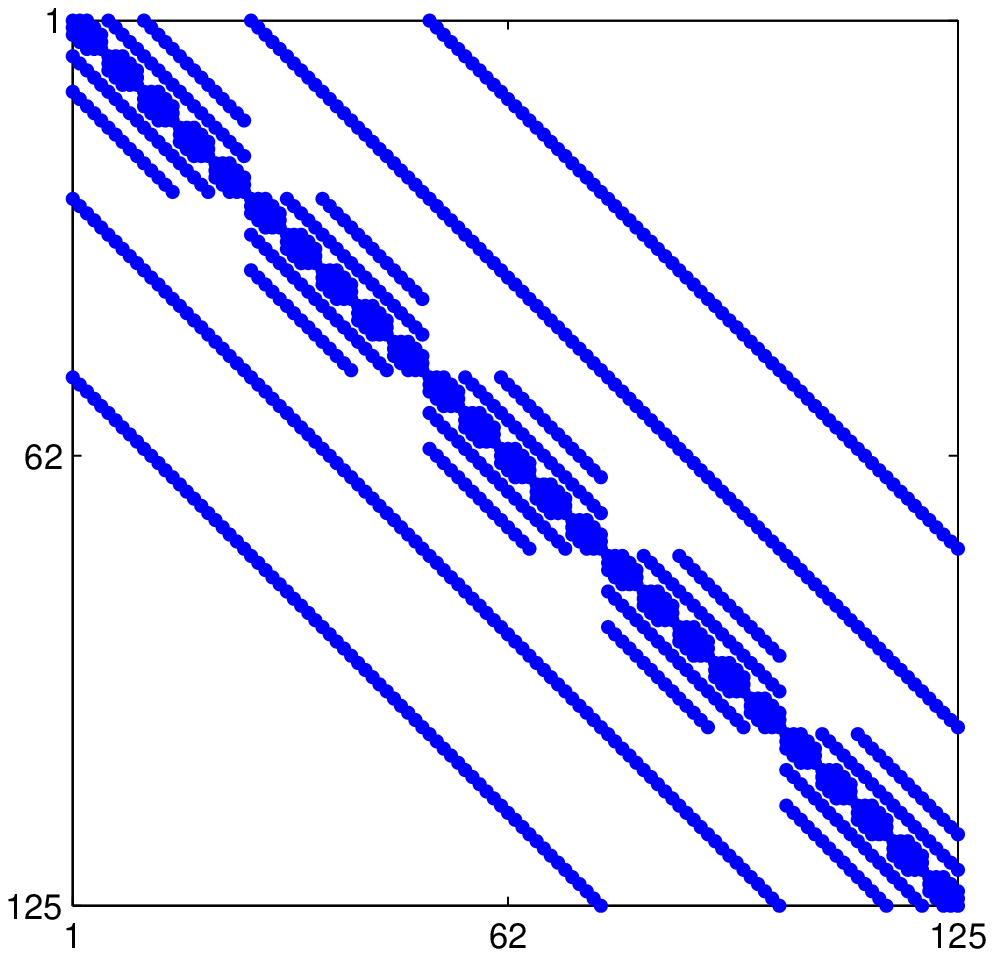}
\caption{(Color online) Form of scheme matrix $A$ and $A^{\prime}$ for the fourth-order 2SHOC scheme of the three-dimensional NLSE.  The dots represent non-zero entries of the matrices.  The matrices shown are for a $7 \times 7$ grid. \label{f:A3D2SHOC}}
\end{figure}

The fourth-order central difference scheme in three dimensions is given by 
\begin{alignat}{2}
\label{3dnonhoc}
\nabla^2 \Psi &\approx \frac{1}{12\,h^2}\left[\Psi_{i+2,j,k} + \Psi_{i-2,j,k} + \Psi_{i,j+2,k} + \Psi_{i,j-2,k} + \Psi_{i,j,k+2} + \Psi_{i,j,k-2}\right. \\
\; &\left. - 16\,(\Psi_{i+1,j,k} + \Psi_{i-1,j,k} + \Psi_{i,j+1,k} + \Psi_{i,j-1,k} + \Psi_{i,j,k+1} + \Psi_{i,j,k-1}) + 90\,\Psi_{i,j,k}\right]. \notag
\end{alignat}
The single-storage version of the 2SHOC equivalent scheme is defined as \cite{ME_2SHOC}
\begin{alignat}{5}
&1) \; D_{i,j,k} = \label{3d2shocs}
\\
&\begin{tabular}{l} 
$\dfrac{1}{h^2}\left(\;
\begin{tabular}{|c|c|c|} \hline
  &   &   \\ \hline
 ~ & 1 & ~  \\ \hline
  &   &   \\ \hline
\end{tabular}
\;\Psi_{i,j+1,k} +
\begin{tabular}{|c|c|c|} \hline
  &  1 &   \\ \hline
1 & -6 & 1 \\ \hline
  &  1 &   \\ \hline
\end{tabular}
\;\Psi_{i,j,k} +
\begin{tabular}{|c|c|c|} \hline
  &   &   \\ \hline
 ~ & 1 & ~  \\ \hline
  &   &   \\ \hline
\end{tabular}
\;\Psi_{i,j-1,k}\;\right),$
\end{tabular} \notag
\\
\;&\; \notag
\\
&2) \; \nabla^2\Psi_{i,j,k} \approx \label{3d2shocs2} 
\\
&\begin{tabular}{l}
$-\dfrac{1}{12}\left(\;
\begin{tabular}{|c|c|c|} \hline
  &   &   \\ \hline
 ~ & 1 & ~ \\ \hline
  &   &   \\ \hline
\end{tabular}
\;D_{i,j+1,k} +
\begin{tabular}{|c|c|c|} \hline
  &   1 &   \\ \hline
1 & -10 & 1 \\ \hline
  &   1 &   \\ \hline
\end{tabular}
\;D_{i,j,k} +
\begin{tabular}{|c|c|c|} \hline
  &   &   \\ \hline
~  & 1 & ~  \\ \hline
  &   &   \\ \hline
\end{tabular}
\;D_{i,j-1,k} \; \right)$
\end{tabular} \notag
\\
&\begin{tabular}{l}
$+ \dfrac{1}{6\,h^2}\left(\;
\begin{tabular}{|c|c|c|} \hline
  & 1 &   \\ \hline
1 &   & 1 \\ \hline
  & 1 &   \\ \hline
\end{tabular} 
\;\Psi_{i,j+1,k} +
\begin{tabular}{|c|c|c|} \hline
1 &     & 1  \\ \hline
  & -12 &    \\ \hline
1 &     & 1  \\ \hline
\end{tabular}
\;\Psi_{i,j,k} +
\begin{tabular}{|c|c|c|} \hline
  & 1 &   \\ \hline
1 &   & 1 \\ \hline
  & 1 &   \\ \hline
\end{tabular}
\;\Psi_{i,j-1,k} \; \right),$
\end{tabular} \notag
\end{alignat}
and the structure of the corresponding $A$ and $A^{\prime}$ matrix are given in Fig.~\ref{f:A3D2SHOC}, while the unique forms of the Gershgorin disk centers and radii for $A^{\prime}$ are shown in Table~\ref{t:3d2shocgd}.
\begin{table}[htbp] 
\centering 
\caption{Gershgorin disk centers and radii for the $A^{\prime}$ matrix of the three-dimensional 2SHOC scheme.}
\begin{center}
\begin{tabular}{|l|l|} \hline
$a_{ii}$    & $r_i = \sum_{i\ne j} |a_{ij}|$ \\ \hline
$L_i - 15/2$   & $33/4$ \\
$L_i - 15/2$   & $25/3$ \\
$L_i - 15/2$   & $101/12$ \\
$L_i - 15/2$   & $17/2$ \\
$L_i - 22/3$   & $67/12$ \\
$L_i - 22/3$   & $17/3$ \\
$L_i - 29/4$   & $17/4$ \\
$L_i - 89/12$  & $83/12$ \\
$L_i - 89/12$  & $7$ \\
$L_i - 89/12$  & $85/12$ \\
\hline
\end{tabular}
\end{center}
\label{t:3d2shocgd}
\end{table}
The stability bounds are the same as in Eq.~(\ref{kbg}), but with $\vec G$ now defined as
\begin{equation}
\label{G3D2SHOC}
\vec G = \frac{1}{12} \times \left\{192,191,190,189,174,173,172,156,155,138,36,21,20,6,5,4,-9,-10,-11,-12 \right\}.
\end{equation}
In the linear case ($s=0$) with no external potential and periodic, Dirichlet or Laplacian-zero boundary conditions, we get the linear stability bound
\begin{equation}
\label{kbound3d2shoclin}
{\bf k} < \frac{\sqrt{8}}{16}\,\frac{h^2}{a},
\end{equation}
which, as in the one- and two-dimensional cases, is simply $3/4$ of the second-order bound given in Eq.~(\ref{kbound3dcdlin}).

\section{Numerical examples}
\label{s:num}
Here we show some numerical examples to demonstrate the accuracy of the predicted stability bounds.  Since the accuracy of the bounds are highly dependent on the problem including the values of $s$, $a$, and $V({\bf r})$, the tests given here are not exhaustive, but serve as a good indication of the general accuracy of the bounds.  

\begin{figure}[p]
 \includegraphics[width=2in]{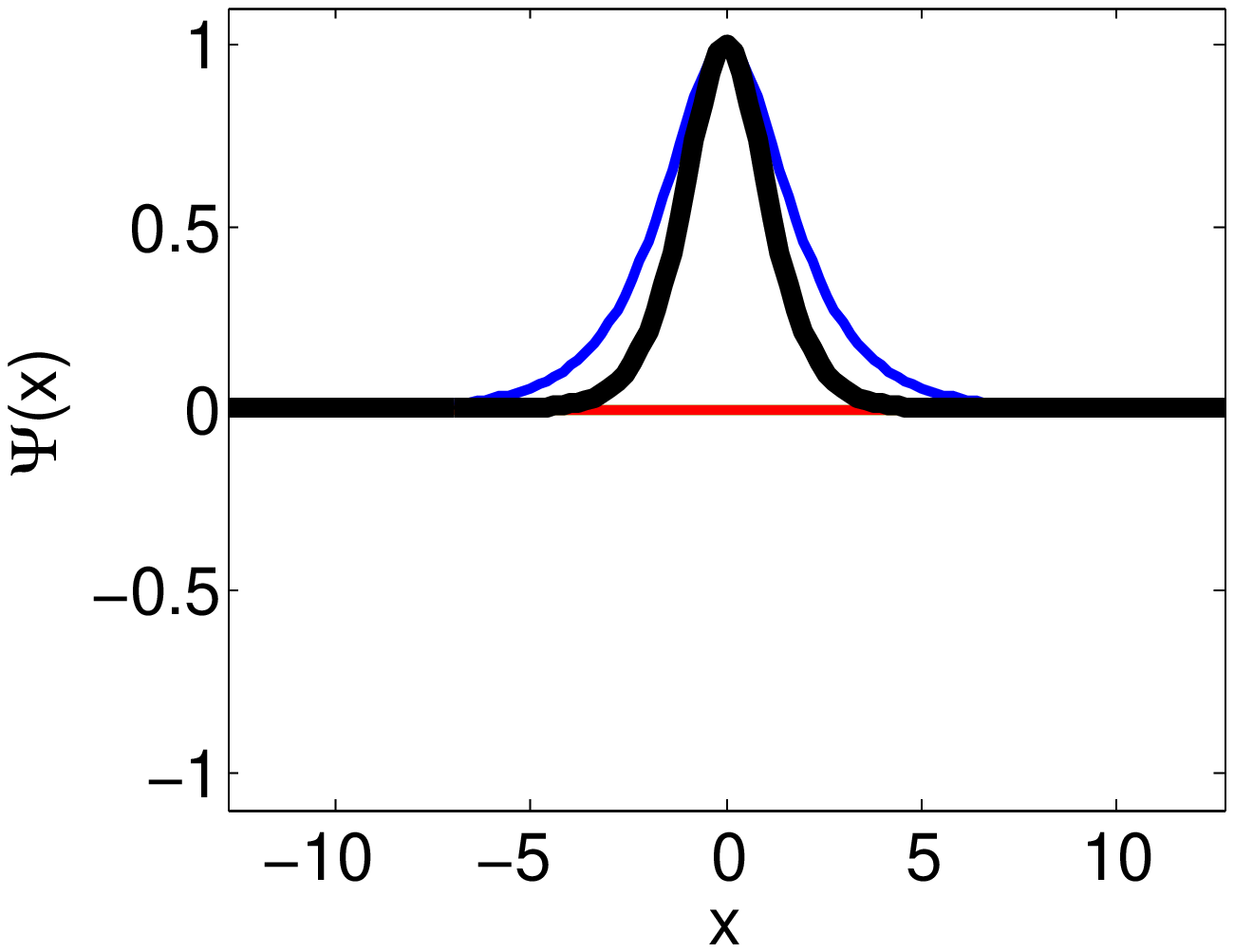}
 \includegraphics[width=2in]{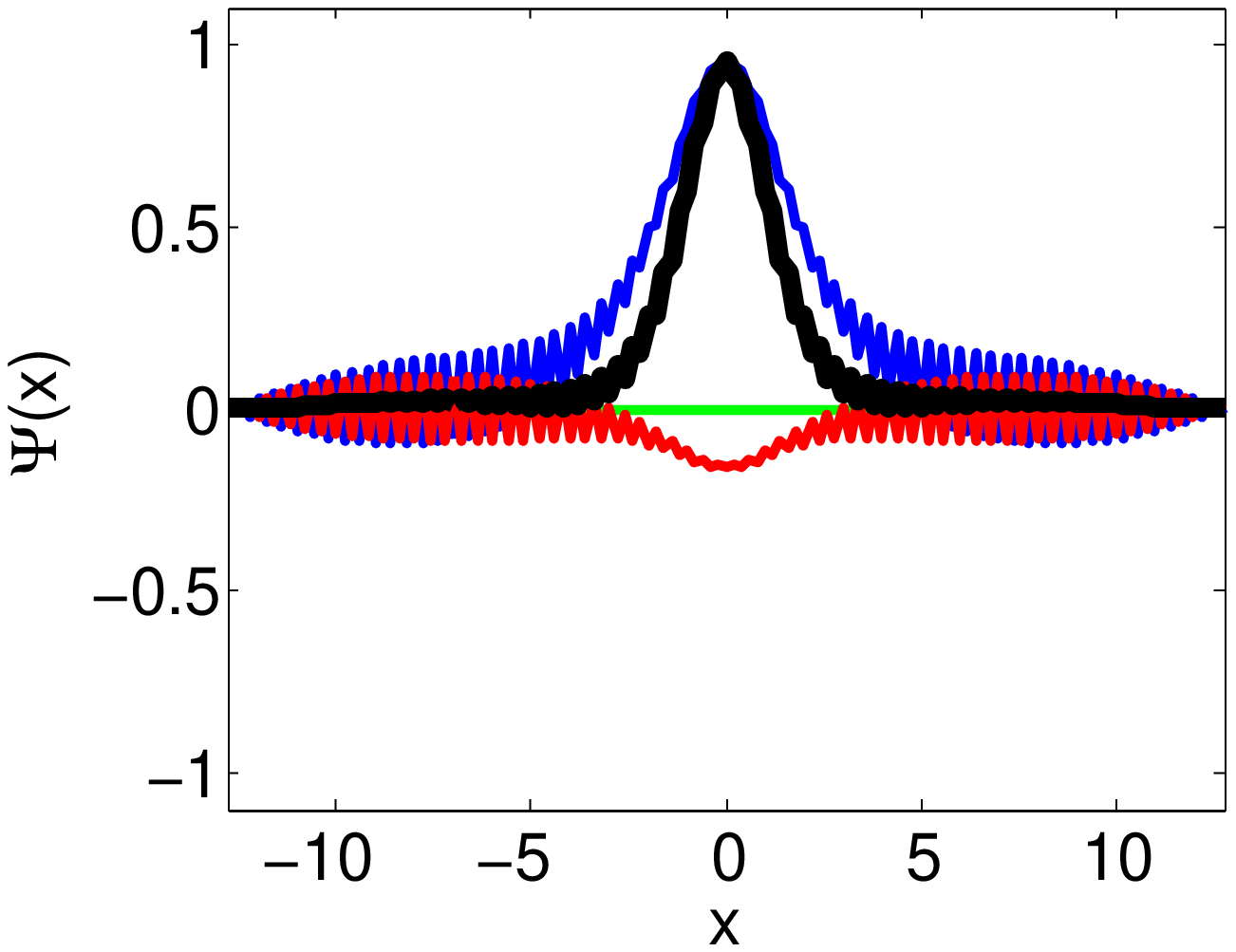}
 \includegraphics[width=2in]{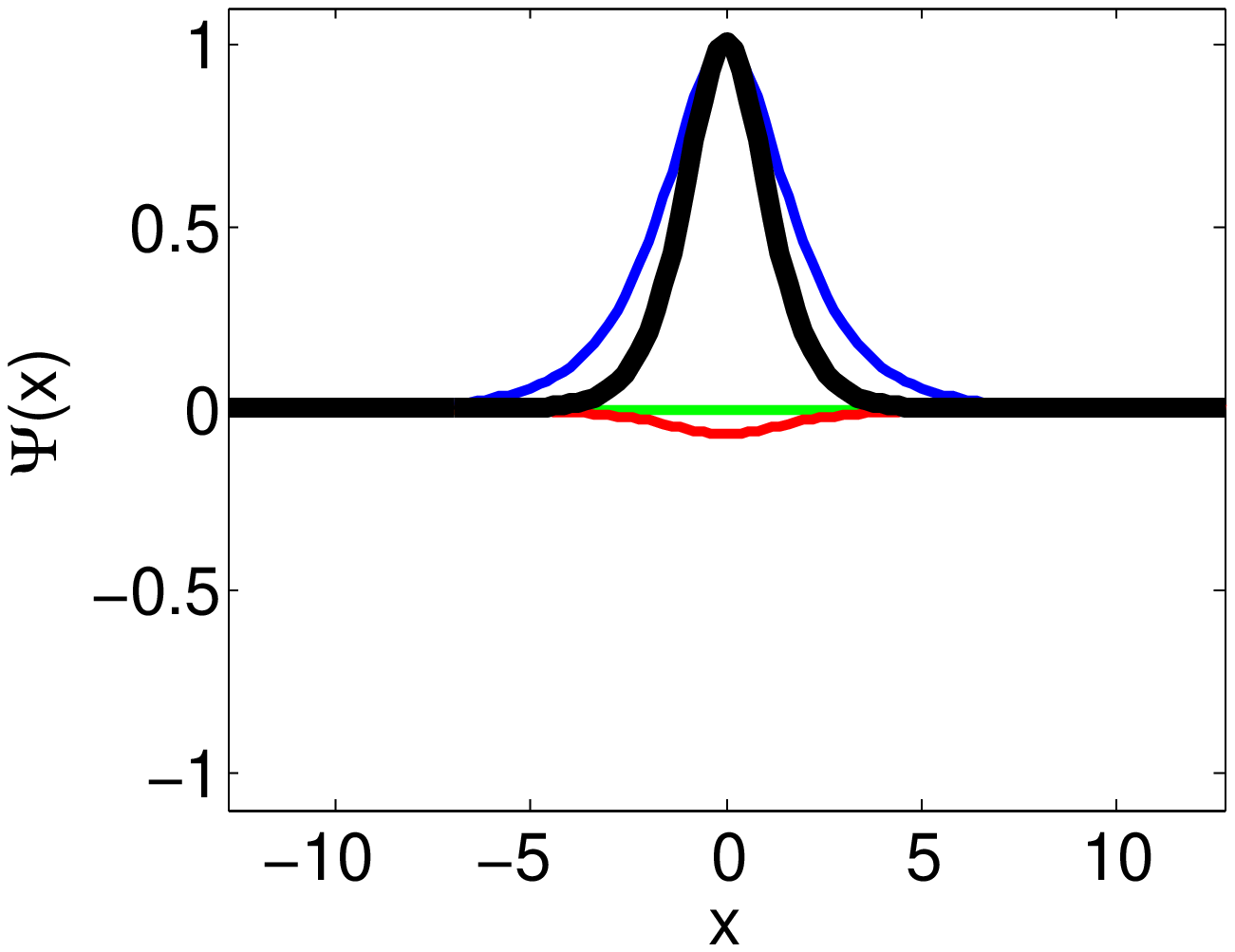}
 \\
 \includegraphics[width=2in]{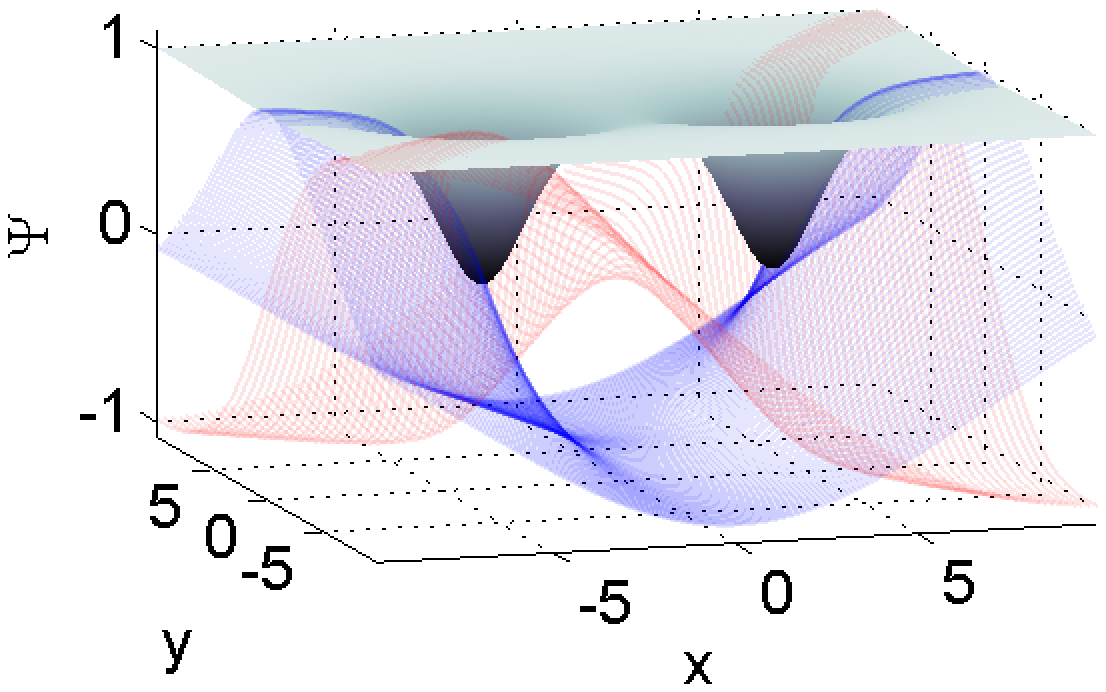}
 \includegraphics[width=2in]{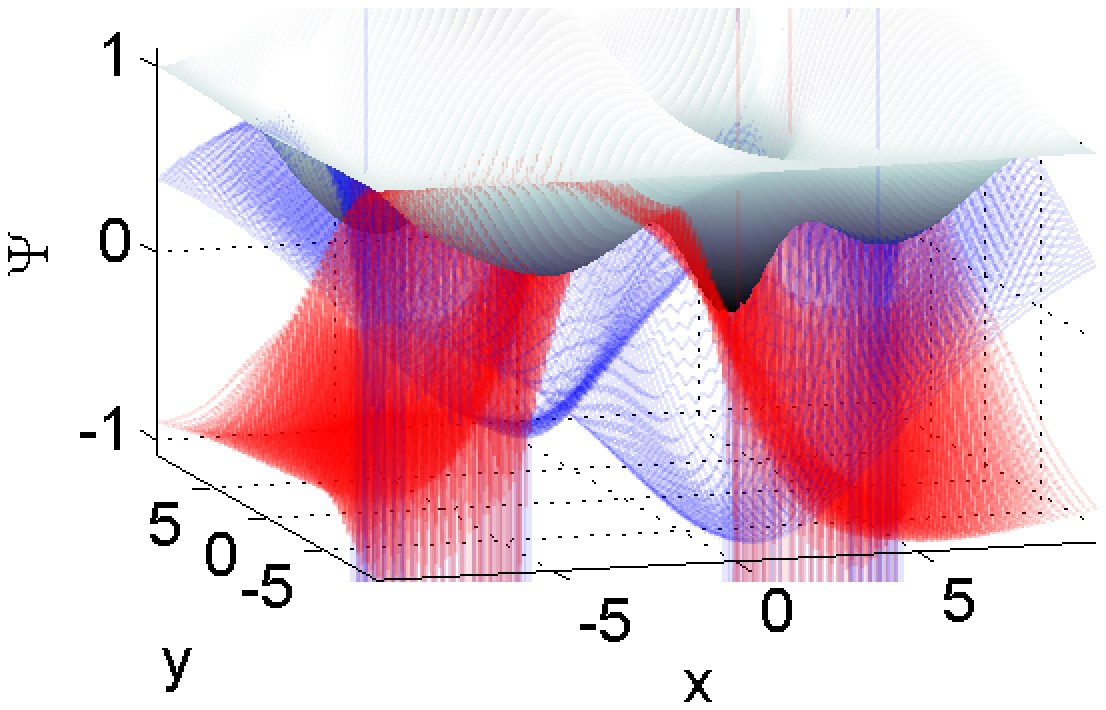}
 \includegraphics[width=2in]{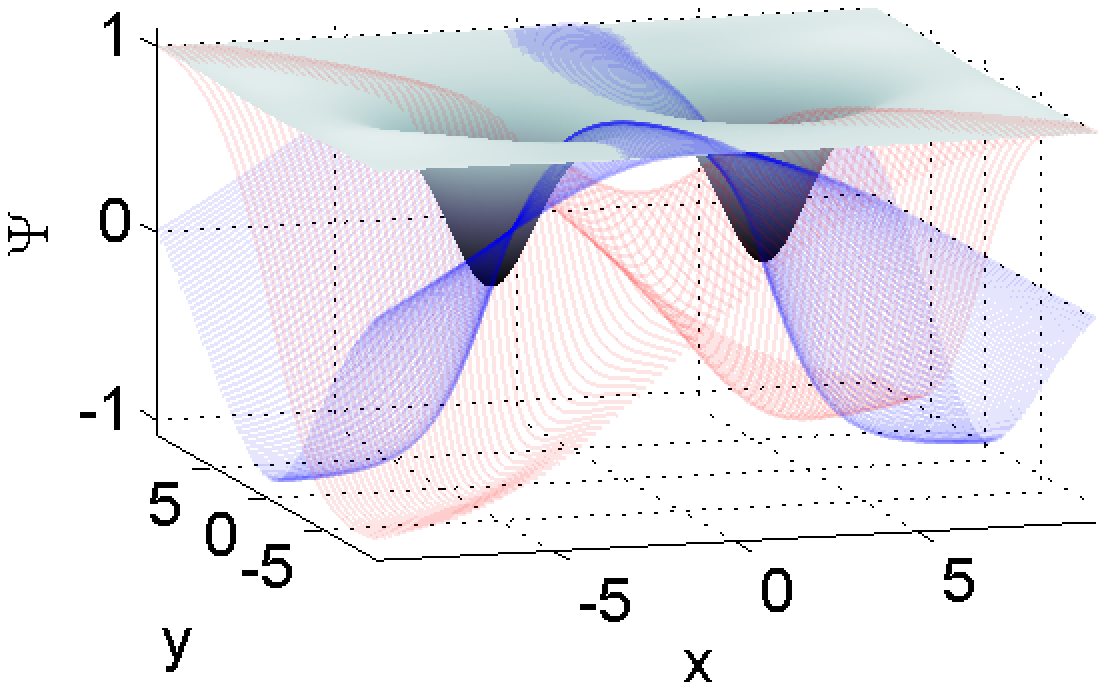}
 \\[2.0ex]
 \includegraphics[width=2in]{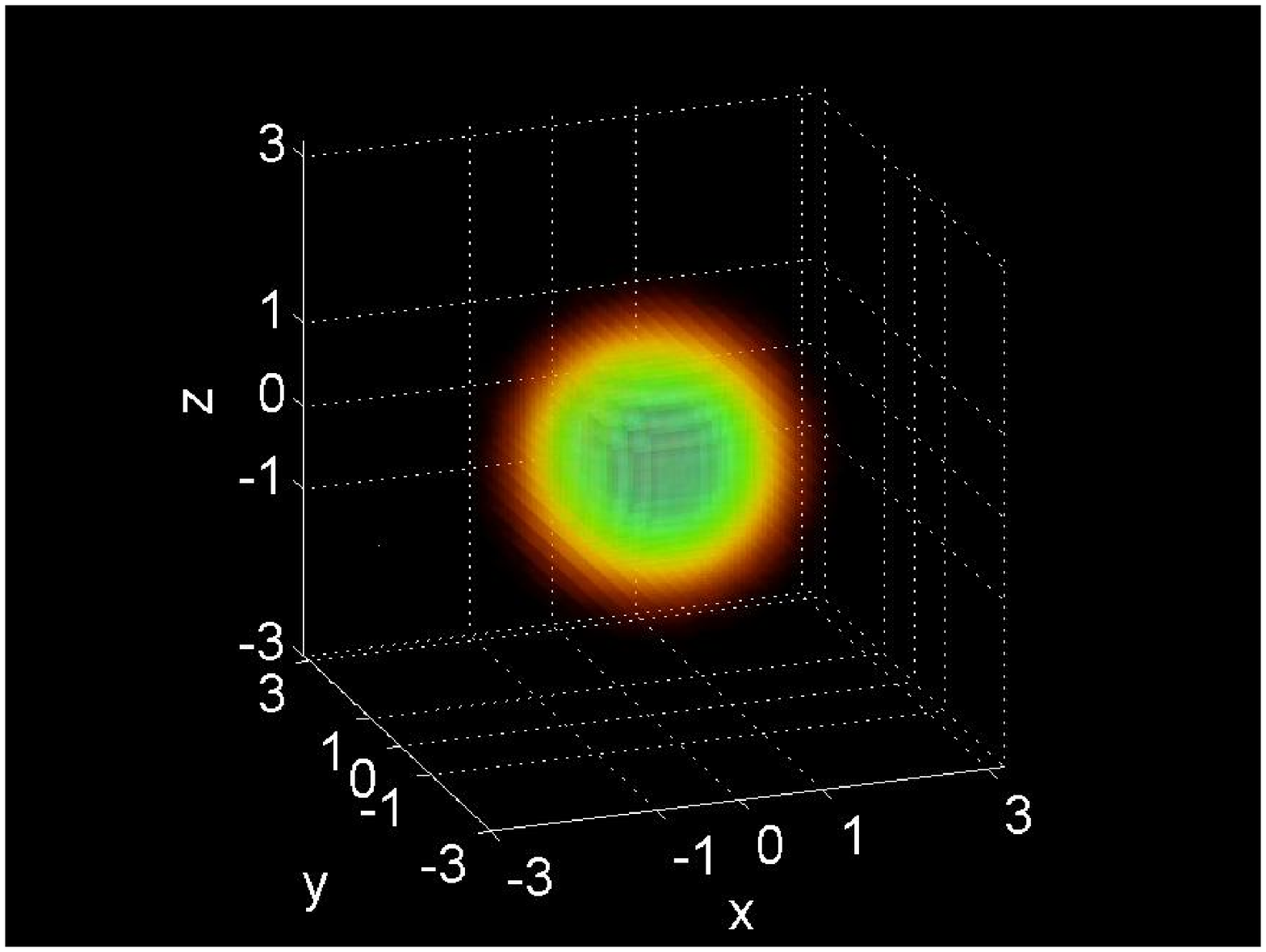}
 \includegraphics[width=2in]{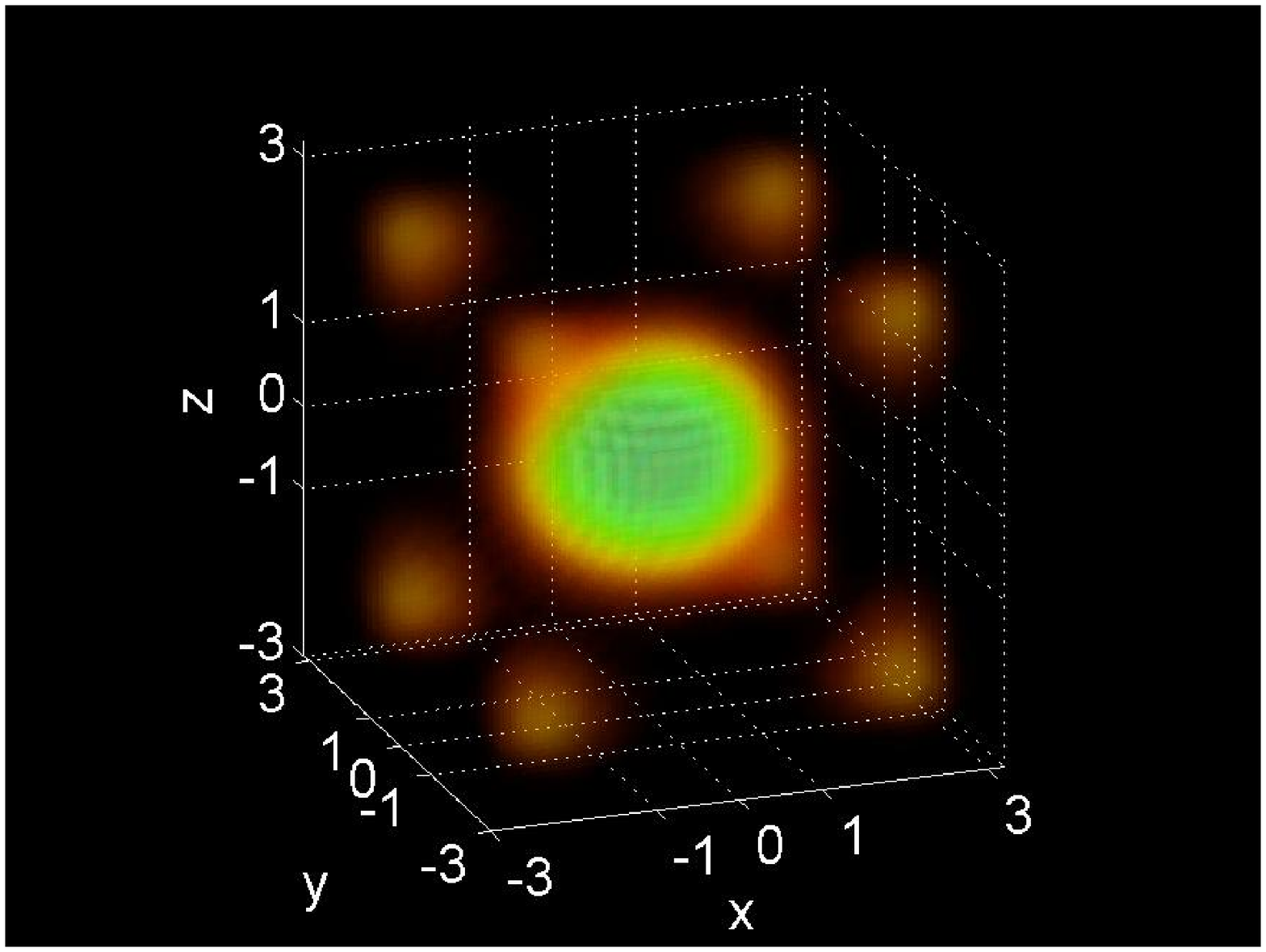}
 \includegraphics[width=2in]{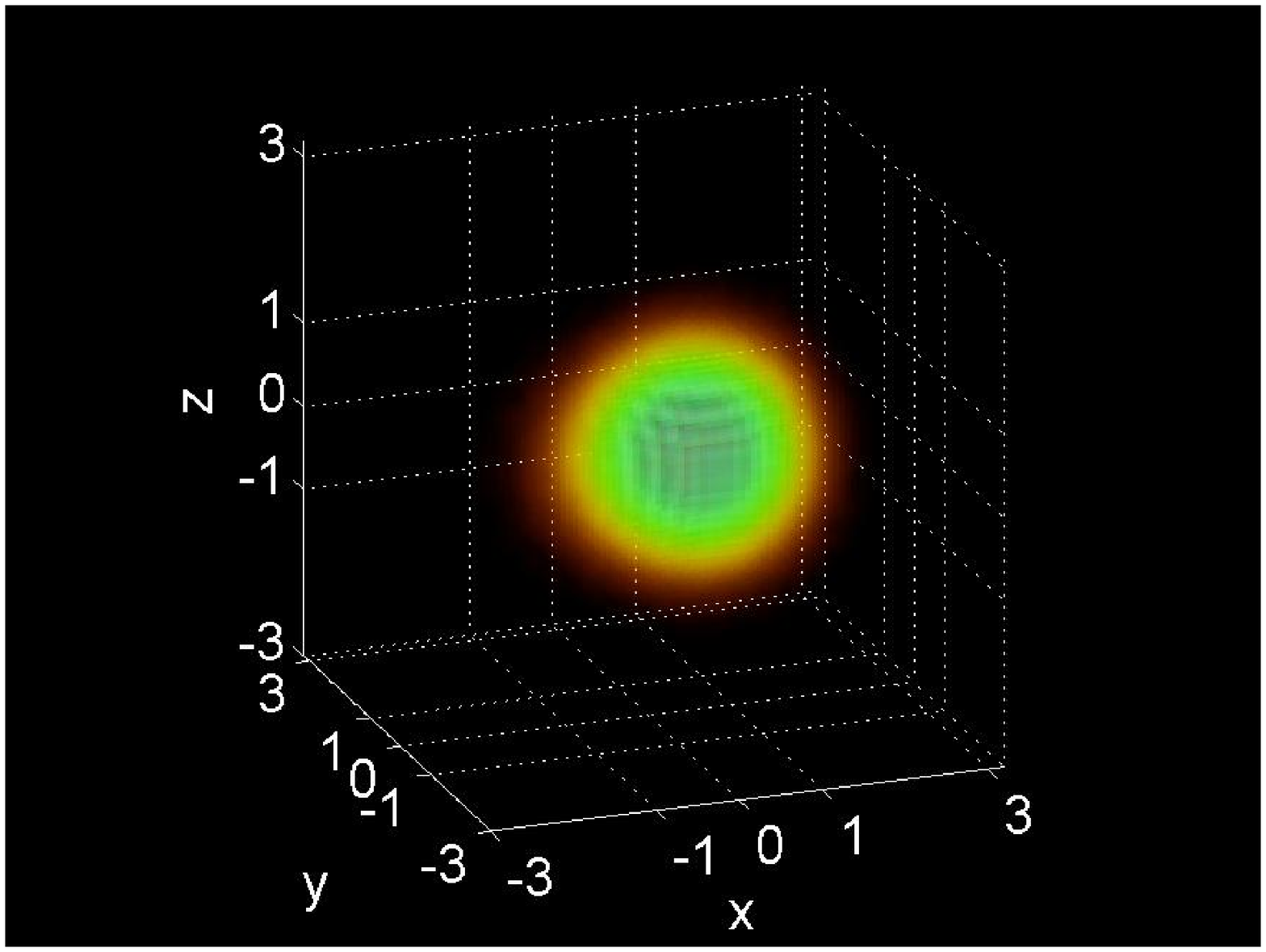}
 \caption{(Color online) Examples of integrating the LSE and NLSE before and after the numerical stability bound for the examples described in Sec.~\ref{s:num}.  Left to right columns:  Initial condition, solution with $k > k_{\mbox{\scriptsize num}}$, and solution with $k = k_{\mbox{\scriptsize num}}$.  Top to bottom:  one-, two-, and three-dimensional test cases using the CD, 2SHOC, and CD schemes respectively. For the one-dimensional test, the predicted stability bounds are $k_{\mbox{\scriptsize lin}}=k_{\mbox{\scriptsize linz}}=0.02828$ and the solution is shown with $k=0.02833$ (middle) and $k=0.02832$ (right) at $t=100$. 
For the two-dimensional test, the predicted stability bounds are $k_{\mbox{\scriptsize lin}}=0.01061$ and $k_{\mbox{\scriptsize linz}}=0.01057$.  The solution is shown with $k=0.01055$ (middle) and $k=0.01054$ (right) at $t=30$ and 
$t=100$ respectively. 
For the three-dimensional test, the predicted stability bounds are $k_{\mbox{\scriptsize lin}}=0.009428$ and $k_{\mbox{\scriptsize linz}}=0.008650$.  The solution is shown with $k=0.009214$ (middle) and $k=0.009213$ (right) at $t=100$.\label{f:results}}
\end{figure}

We choose to use three initial conditions, one for each dimensionality case,
and integrate them using both the CD and 2SHOC schemes.  
In one dimension, we use the exact steady-state bright soliton solution to the NLSE with $V(x)=0$ and $s>0$ given as \cite{SOL_Bright_Gray_Dark_Opt}
\begin{equation}
\label{soliton}
\Psi(x,t) = \sqrt{\frac{2\,\Omega}{s}}\,\mbox{sech}\left[\sqrt{\frac{\Omega}{a}}\,x\right]\,\mbox{exp}\left(i\,\Omega\,t\right),
\end{equation} 
where $\Omega$ is the frequency, and we set $V(x)=0$, $\Omega=1$, $s=1$, and $a=1$ and use Dirichlet boundary conditions ($\Psi = 0$).  In two dimensions, we use an approximation to a co-rotating dark vortex pair solution to the NLSE.  Each vortex is given by \cite{BEC_DYN_NONLIN}
\begin{equation}
\label{exmp2Dvort}
\Psi(r,\theta,t) = f(r)\,\mbox{exp}[i\,(m\,\theta + \Omega\,t)],
\end{equation}
where $m$ is the topological charge of the vortex (in our case, we use $m=1$), $\Omega = -1$, and we use $s=1$ and $a=1$.  The term $f(r)$ is the real-valued radial profile centered at the vortex core which can be found numerically \cite{ME_DISS}.  The pair of vortices are then combined to yield the initial condition
\[
\Psi(x,y,t=0) = f(r_1)f(r_2)\,\mbox{exp}[i\,m\,(\theta_1 + \theta_2)],
\]
where $r_1=\sqrt{(x-x_0)^2 + y^2}$, $r_2=\sqrt{(x+x_0)^2 + y^2}$, $\theta_1=\tan^{-1}(y/(x-x_0))$, and $\theta_2=\tan^{-1}(y/(x+x_0))$, which approximates the true initial condition of a co-rotating steady-state vortex pair solution
located at $(-x_0,0)$ and $(x_0,0)$. Here we choose $x_0=4$.  
Since $|\Psi|^2$ does not decay at the boundaries, we use the modulus-squared Dirichlet boundary condition $|\Psi|^2 = 1$.
In three dimensions, we use a steady-state bright Gaussian solution of the LSE in a potential trap with an added initial `kick' in the $x$-direction which causes the structure to oscillate in the $x$-direction.  The initial condition is given by
\begin{equation}
\label{3dexpsol}
\Psi(x,y,z,t=0) = \mbox{exp}\left(-\frac{x^2 + y^2 + z^2}{2\,a}\right) \mbox{exp}\left(-i\,\frac{x}{2}\right),
\end{equation}
with external potential
\begin{equation}
\label{3dexpsolv}
V(x,y,z) = \frac{x^2 + y^2 +z^2}{a},
\end{equation}
and we use the Dirichlet boundary condition $\Psi = 0$.  
For all simulations, we set the spatial step-size of the grid to $h=1/5$.  The three initial conditions are shown in the left column of Fig.~\ref{f:results}.

To test the stability bounds, each solution is integrated to an ending time of $t=100$ and it is observed if the solution remains stable.  We increase the time-step $k$ until the solution becomes unstable within the $t=100$ simulation time, at which point the largest time-step that was stable is denoted $k_{\mbox{\scriptsize num}}$.  This is then compared to the computed linear [Eqs.~(\ref{kdcdlin}) and (\ref{kd2shoclin}) denoted $k_{\mbox{\scriptsize lin}}$] and fully linearized [Eq.~(\ref{kdfull}) denoted $k_{\mbox{\scriptsize linz}}$] stability bounds formulated in Secs.~\ref{s:1dstb}--\ref{s:3dstb}.  The time-step is incremented to yield the numerical stability limit to within four significant figures.  All of the simulations are performed using the NLSEmagic software package \cite{NLSEMAGIC}.
  
\begin{table}[p]
\begin{center}
\begin{tabular}{|l|r|r|r|r|r|} \hline
Example: & $k_{\mbox{\scriptsize lin}}$ & $k_{\mbox{\scriptsize linz}}$ & $k_{\mbox{\scriptsize num}}$ & \%-diff $k_{\mbox{\scriptsize lin}}$ & \%-diff $k_{\mbox{\scriptsize linz}}$ \\
\hline
1D CD    & 0.02828  & 0.02828  & 0.02832  &  0.14 &  0.14 \\
1D 2SHOC & 0.02121  & 0.02121  & 0.02124  &  0.14 &  0.14 \\
2D CD    & 0.01414  & 0.01407  & 0.01402  & -0.85 & -0.36 \\
2D 2SHOC & 0.01061  & 0.01057  & 0.01054  & -0.66 & -0.28 \\
3D CD    & 0.009428 & 0.008650 & 0.009213 & -2.28 &  6.51 \\
3D 2SHOC & 0.007071 & 0.006624 & 0.006992 & -1.12 &  5.56 \\
\hline
\end{tabular}
\caption{Numerical test results of finding the numerical stability bound ($k_{\mbox{\scriptsize num}}$) for the example problems described in Sec.~\ref{s:num} compared to the predicted linear ($k_{\mbox{\scriptsize lin}}$) and linearized ($k_{\mbox{\scriptsize linz}}$) bounds.\label{t:results}}
\end{center}
\end{table}

Before displaying the results, we point out that there are some sources of error to consider.  First, the predicted stability bounds are linearized and therefore will not be the same as the corresponding true nonlinear stability bounds.  Second, in our analysis, we chose to use every possible combination of $\vec L$ which may lead to predictions of the bounds which are stricter than the true bound.  Finally, it is sometimes difficult to determine the true stability bound numerically, as some unstable time-steps may only exhibit their instability after a very long simulation time.  For our test, we choose a moderately long simulation time, but the exact bound may be slightly higher than the given result. 

The results are shown in Table~\ref{t:results}, while Fig.~\ref{f:results} shows the solutions before and after the recorded numerical stability bounds for three chosen examples.
We see that overall, the numerical results match the predicted stability 
values quite well (especially in one and two dimensions) with a maximum 
percent difference of $6.5\%$ when $V({\bf r})\not= 0$ in the three-dimensional
example, but with a typical percent difference less that $1\%$ when 
$V({\bf r})=0$.  
It is noted that in some cases the predicted bounds are stricter than the numerical result, while in other cases, they are too lenient, noting that the examples with $s>0$ were all too strict, while those with $s<0$ were all too lenient.  However, due to the small number of tests, no conclusions about the effect of the sign and presence of the parameters and external potential of the LSE and NLSE on the stability bound predictions can be drawn from these observations.  
In terms of choosing a stable time-step for LSE and NLSE simulations, the results given are well within a tolerable range, and in practice one would use a time-step some percentage (say $10$--$20\%$) lower than the predicted bound to ensure stability over long integration times.

\section{Conclusion and summary of results}
\label{s:sum}
In this paper we have formulated linearized stability bounds for using second- and fourth-order spatial finite-differencing with fourth-order Runge-Kutta time-stepping for the multi-dimensional nonlinear Schr{\"o}dinger equation (NLSE) with Dirichlet, modulus-squared Dirichlet, Laplacian-zero, and periodic boundary conditions.  

A summary of the stability results for easy reference is given presently.  For the nonlinear Schr{\"o}dinger equation defined as 
\[
i\, \frac{\partial \Psi}{\partial t} +a\,\nabla^2 \Psi -V({\bf r})\Psi + s\,|\Psi|^2\Psi = 0,
\]
where $a>0$ and $s$ are parameters of the system and $V({\bf r})$ is an external potential, the numerical stability bounds on the time-step when using the fourth-order Runge-Kutta time-stepping scheme is as follows: 

In the linear case where $s=0$ and with no external potential ($V({\bf r})=0$), utilizing periodic, Dirichlet, or Laplacian-zero boundary conditions, the stability bound on the time-step ${\bf k}$ when using the second-order central difference (CD) scheme in a $d$-dimensional setting is 
\begin{equation}
\label{kdcdlin}
{\bf k_{\mbox{\scriptsize CD}}} < \frac{h^2}{d\,\sqrt{2}\, a},
\end{equation}
while that of using a fourth-order central difference scheme (with interior points computed in the two-step high-order compact (2SHOC) methodology of Ref.~\cite{ME_2SHOC}) is
\begin{equation}
\label{kd2shoclin}
{\bf k_{\mbox{\scriptsize 2SHOC}}} < \left(\frac{3}{4}\right) \frac{h^2}{d\,\sqrt{2}\, a}.
\end{equation}
The linearized stability bounds for the general NLSE are
\begin{equation}
\label{kdfull}
{\bf k} < \frac{\sqrt{8}}{\max\{\lVert\vec B\rVert_{\infty},\Vert \forall L_i, L_i-\vec G\rVert_{\infty}\}}\,\frac{h^2}{a},
\end{equation}
where $\vec B$ are the boundary points as defined by Table~\ref{t:bc2} (or in the periodic case is ignored), the elements of $\vec L$ is defined as
\[
L_i = \frac{h^2}{a}\left(s|\Psi_i|^2 - V_i\right),
\]
where the index $i$ spans the entire grid, and $\vec G$ is a set of values defined in Table~\ref{t:sumresults}, determined by the dimension and method being used.
\begin{table}[htbp] 
\caption{Values of $\vec B$ in Eq.~(\ref{kdfull}).}
\begin{center}
\begin{tabular}{|l|c|c|c|} \hline
$\;$  & Dirichlet ($\Psi_b = \mbox{const}$) & Laplacian-zero ($\nabla^2\Psi_b = 0$)        & MSD ($|\Psi_b|^2 = \mbox{const}$) \\ \hline
$B_b$ & $0$       & $\dfrac{h^2}{a}\left(s|\Psi_b|^2-V_b\right)$ & $\dfrac{h^2}{a}\mbox{Im}\left[\dfrac{\Psi_{t,b-1}}{\Psi_{b-1}}\right]$ \\ \hline
\end{tabular}
\end{center}
\label{t:bc2}
\end{table}
\begin{table}[htbp] 
\caption{Values of $\vec G$ in Eq.~(\ref{kdfull}).}
\begin{center}
\begin{tabular}{|l|c|c|} \hline
Scheme $\rightarrow$ & CD $O(h^2)$         & 2SHOC $O(h^4)$ 
\\ \hline
\;     & \;                  &\;              
\\
1D     & $\{4,3,1,0\}$       &$\dfrac{1}{12} \times \left\{64,63,46, \right.$     
\\
       & \;                  &$\left. \qquad 12,-3,-4 \right\}$
\\
2D     & $\{8,7,6,2,1,0\}$   &  $\dfrac{1}{12} \times \left\{128,127,126,110,109,\right.$ 
\\
\;     &  \;                 & $\left.  \qquad 92,24,9,8,-6,-7,-8 \right\} $
\\
3D     & $\{12,11,10,9,3,2,1,0\}$  & $\dfrac{1}{12} \times \left\{192,191,190,189,174,\right.$
\\
\;     & \;                        & $\qquad 173,172,156,155,138,36,21,$
\\
\;     & \;                        & $\qquad \left. 20,6,5,4,-9,-10,-11,-12 \right\}$  
\\ \hline
\end{tabular}
\end{center}
\label{t:sumresults}
\end{table}

We have found through numerical testing (those of Sec.~\ref{s:num}, as well as others not reported here) that to ensure stability in all dimensions for typical problems, the bounds must be lowered by about $10\%$--$20\%$ (most likely due to nonlinear effects).  Also, we note that the reduced linear results are often similar to the full linearized bounds and can therefore be used as a good quick estimate of the stability bound.

\section*{Acknowledgments}
This research was supported by NSF-DMS-0806762 and the Computational
Science Research Center (CSRC) at SDSU.
We gratefully acknowledge insightful discussions with Peter Blomgren.

\def\myitemsep{5pt}
\bibliographystyle{elsart-num-sort}
\bibliography{RK4STB5}  
\end{document}